\documentclass[12pt]{article}
\usepackage[a4paper,margin=2cm]{geometry}

\usepackage{xr} 
\usepackage{tikz}
\usetikzlibrary{decorations.pathreplacing,calligraphy}
\usetikzlibrary{decorations.markings}
\usetikzlibrary{patterns}
\usepackage{xfrac}

\usepackage{makecell}
\usepackage[utf8]{inputenc}
\usepackage[english]{babel}
\usepackage{amssymb}
\usepackage{amsmath}
\usepackage{graphicx}
\usepackage{caption} 
\usepackage{dsfont}
\usepackage{fancyhdr}
\usepackage{amsthm}
\usepackage[T1]{fontenc}
\usepackage{fancyhdr}
\usepackage{hyperref}
\usepackage{array}

\usepackage{minitoc}
\usepackage{geometry} 
\usepackage{nomencl}
\usepackage{url}
\usepackage{subcaption}
\usepackage{wallpaper}

\usepackage{cancel}
\usepackage{stmaryrd}

\usepackage{float}
\usepackage{geometry} 

\usepackage[shortlabels]{enumitem}

\graphicspath{ {figures/} }
\usepackage{subcaption}
\usepackage{wallpaper}

\newtheorem{theorem}{Theorem}[section]

\newtheorem{lemma}[theorem]{Lemma}
\newtheorem{prop}[theorem]{Proposition}

\theoremstyle{definition}

\theoremstyle{remark}
\newtheorem{remark}[theorem]{Remark}

\newcommand{\vd}{\,\mathrm{d}}

\newcommand{\Rp}{(0,\infty)}
\newcommand{\Lk}{\mathcal{L}}
\newcommand{\qLk}{q\text{-}\mathcal{L}}
\newcommand{\calH}{\mathcal{H}}
\newcommand{\calM}{M}
\newcommand{\frM}{\mathcal{M}}
\newcommand{\calQ}{Q}
\newcommand{\frQ}{\mathbf{Q}}
\newcommand{\Emu}{\mathbb{E}}
\newcommand{\Xn}{X^{(n)}}

\usepackage[numbers]{natbib}

\usepackage{amsmath}
\numberwithin{equation}{section}

\usepackage{hyperref}
\usepackage{cleveref}
\usepackage{autonum}

\title{Parameters estimation of a \\ Threshold Chan-Karolyi-Longstaff-Sanders process \\ from continuous and discrete observations}
\author{  
Sara Mazzonetto\thanks{Université de Lorraine, CNRS, Inria, IECL, F-54000 Nancy, France, sara.mazzonetto@univ-lorraine.fr}
        \and
        Benoît Nieto\thanks{ Centre de Mathématiques Appliquées (CMAP), CNRS, Ecole Polytechnique, Institut Polytechnique de Paris, Palaiseau, France,
benoit.nieto@polytechnique.edu}         }

\begin{document}
\maketitle    
\noindent \textbf{Abstract.}  We consider a continuous time process that is self-exciting and ergodic, called threshold Chan-Karolyi-Longstaff-Sanders (CKLS) process. This process is a generalization of various models in econometrics, such as Vasicek model, Cox-Ingersoll-Ross, and Black-Scholes, allowing for the presence of several thresholds which determine changes in the dynamics. We study the asymptotic behavior of maximum-likelihood and quasi-maximum-likelihood estimators of the drift parameters in the case of continuous time and discrete time observations. We show that for high frequency observations and infinite horizon the estimators satisfy the same asymptotic normality property as in the case of continuous time observations. We also discuss diffusion coefficient estimation. 
Finally, we apply our estimators to simulated and real data to motivate considering (multiple) thresholds.


\medskip


\section{Introduction}
We consider a stochastic process that follows different Chan-Karolyi-Longstaff-Sanders (CKLS) dynamics on separate intervals. It satisfies the following stochastic differential equation (SDE):
\begin{equation}
X_t = X_0 + \int_0^t (a(X_s) - b(X_s)X_s) \vd s + \int_0^t \sigma(X_s)|X_s|^{\gamma(X_s)} \vd B_s, \quad t \geq 0,
\end{equation}
where the coefficients $a$, $b$, $\sigma$, $\gamma$ are piecewise constant functions which are defined formally in Section \ref{Preliminaries}. We refer to this model as the \emph{threshold CKLS} (T-CKLS) process. The \emph{thresholds} are the discontinuity points of the coefficients. 

The class of T-CKLS models includes several celebrated non-threshold models. Before to mention them, let us recall the equation satisfied by (non-threshold) CKLS:
\begin{equation}
    X_t=x_0+\int_0^t (a-b X_s)\vd s+\int_0^t \sigma |X_s|^\gamma \vd B_s
    \label{CKLS}
\end{equation}
where $x_0,a,\sigma\in (0,\infty)$, $b\in \mathbb{R}$, $\gamma\in[0,1]$. CKLS generalizes Ornstein-Uhlenbeck (OU) or Vasicek when $\gamma=0$, Cox-Ingersoll-Ross (CIR) when $\gamma=1/2$, Black-Scholes when $\gamma=1$ and $a=0$, Merton when $\gamma=0$ and $b=0$, Constant Elasticity of Variance (CEV) model when $a=0$. 
These processes are well-known models for instance for interest rate modeling (see, e.g.,~\cite{chan1992empirical} for CKLS). 
Their threshold analogues have also been considered. In particular T-OU process, which corresponds to T-CKLS with $\gamma\equiv 0$, is also called Self-Exciting Threshold (SET) Vasicek model~\cite{decamps2006self} because it excites itself by changing dynamics according to its own position. 
Because of its similarities with threshold autoregressive (TAR) processes, T-OU is also considered as a first-order continuous time TAR process, see~\cite{su2015quasi}. (See, e.g.,~\cite{hansen2011threshold} for a review of TAR processes in economics).

Clearly the T-CKLS and its special cases such as T-OU and T-CIR ($\gamma\equiv 1/2$) are SET processes in the meaning of~\cite{decamps2006self} recalled above.
These processes appear also under the more general nomenclature or \emph{threshold diffusions} (TDs) in~\cite{su2015quasi}, they model regime-switching behaviors in time series econometrics. Indeed TDs, being diffusions with piecewise regular coefficients, allow for 
several mean reversion levels. TDs also attract attention for applications in financial modelling, e.g., in~\cite{decamps2006self, lejay2019pigato, motaesquivel}, population ecology~\cite{cantrell1999diffusion}, etc. 

Recently, several studies have been conducted on the parametric estimation for TDs. Just to mention some results, in the case of continuous time observations~\cite{kutoyants2012, su2015quasi, lejay2020maximum}, high frequency observations on finite or infinite horizon~\cite{lejay2018statistical, mazzonetto2020drift, mazzonetto2024BM}, low frequency ones~\cite{Yu2020, HuXi, mazzonetto2024BM}. 
Non-threshold processes represent a special class of TDs. Let us also mention the literature on parametric estimation for these processes. For CKLS, the drift maximum likelihood estimation has been considered in~\cite{mishura2022parameter} in the continuous time setting, and for CIR let us mention the series of works~\cite{alaya2012parameter, ben2013asymptotic,alaya:hal-03936795} both in the continuous and discrete time setting.

In this paper, we focus on parametric estimation for T-CKLS, which includes most of the TDs considered in the just mentioned literature, and also allows for different dynamics in fixed intervals: for instance, OU in one interval, CIR in another, and possibly other special cases of CKLS in other intervals. 
We provide theoretical results on the estimation of drift and diffusion parameters of \emph{ergodic} T-CKLS assuming that the piecewise constant function $\gamma$ is known, thresholds included. 
We consider two scenarios, one in which we dispose of an entire trajectory of the T-CKLS (continuous time observations) and the more realistic one in which we dispose of the observation at some times (discrete observations). 
Furthermore, we apply our estimators to both simulated and real datasets to assess their practical accuracy. Our numerical study also includes threshold estimation, and we propose a statistical test procedure to detect the presence of one or more thresholds in the drift and volatility of the T-CKLS model. Regarding the estimation of $\gamma$, we assume it is known in our numerical study. For the non-threshold CKLS model, an estimator for $\gamma$ was proposed in~\cite{mishura2022parameter}; however, its properties have not yet been studied, particularly in the context of discrete observations. While this analysis is beyond the scope of our work, it remains necessary for both the CKLS and T-CKLS models, especially in conjunction with threshold estimation. 

Let us describe our results.\\

In the continuous time setting, we assume that the trajectory is observed on a time window $[0,T]$. We consider likelihood and quasi-likelihood functions in~\eqref{Likelihood}-\eqref{Quasi_likelihood} based on the 
Girsanov's weight and compute, in Proposition~\ref{Proposition_estimateur_continu}, the respective maximum likelihood estimator (MLE) and a quasi maximum likelihood estimator (QMLE) in terms of functionals of the trajectory. Note that since T-CKLS follows CKLS dynamics on different intervals, the estimation of the coefficients of the CKLS dynamic in an interval depends only on the part of the trajectory falling in that interval. 
We then study the drift estimators asymptotic properties in \emph{long time}, i.e.~as $T \to +\infty$. Theorem~\ref{Theorem_cont} establishes the consistency (law of large numbers) and asymptotic normality (central limit theorem) of the estimators. 
Its proof is based on classical martingale central limit theorems, as it is the case for other results for ergodic TDs, e.g.~\cite{su2015quasi,lejay2020maximum,mazzonetto2020drift}.
Theorem~\ref{Theorem_cont} recovers these results and generalizes them in a common setting. 
This was already done for QMLE in \cite{su2015quasi} under stronger assumptions, but it is new for MLE.
Theorem~\ref{Theorem_cont} is a key step for the main results of this paper, which are the analogous results for discrete observations.
\\

In the case of discrete observations, we assume that the ergodic process is in its stationary regime (i.e., $X_t$ is distributed according to the stationary distribution for all $t \geq 0$) and it is observed on a time-grid with $N$ observations (not necessarily equally spaced) with maximal lag between two consecutive observations denoted by $\Delta_N$. 
Note that the law of the sample cannot be described to obtain a likelihood function because, except in very special cases of TDs, the law of the process, particularly the finite-dimensional distributions, is not explicitly known. Therefore, we discretize both the likelihood and quasi-likelihood functions from the continuous setting and define our drift estimators as the maximizers of these discretized functions, which we continue to denote as MLE and QMLE.
This is a standard procedure in estimations of diffusions from discrete observations, e.g.,~ \cite{genon1990maximnm,yoshida1992estimation}. Inspired by the MLE estimation of standard CIR process~\cite{ben2013asymptotic}, 
we derive an alternative discretized likelihood using a different discretization method based on Itô-Tanaka formula. This results in a MLE with more easily provable convergence. Additionally, two different discretization methods for the quadratic variation yield an estimator for the diffusion coefficient. 

We study the asymptotic behavior of the above mentioned estimators as the number of observations $N \to +\infty$ under the assumptions that $\lim_{N \to \infty} \Delta_N \to 0$ (high frequency) and $\lim_{N \to \infty} T_N = +\infty$ (long time). 
In this setting, Theorem~\ref{Thm_sigma} and Theorem~\ref{Theorem_discret_cons_tcl} establish the consistency of the diffusion coefficient estimator, QMLE, and MLE. 
A lower bound for the speed of convergence of the diffusion coefficient estimator is obtained in~Theorem~\ref{Thm_sigma} and asymptotic normality of the QMLE and MLE is shown in Theorem~\ref{Theorem_discret_cons_tcl}. The latter results require the additional condition $\lim_{N \to \infty} \Delta_N T_N = 0$. This condition is consistent with the existing literature on parameter estimation for diffusions from discrete observations, with and without threshold effects, e.g.,~\cite{amorino_gloter2020, alaya2012parameter, mazzonetto2020drift, mazzonetto2024BM}. 

The proof strategies of our results is different than the ones employed for MLE estimation for standard SDEs and, in particular, CIR process~\cite{alaya2012parameter, ben2013asymptotic}. Those proofs rely on the explicit knowledge of the law of the process,  
and, as already mentioned, the law of the T-CKLS process is still unknown. 
Our proof strategy is inspired by the one for drift estimation in the T-OU model in~\cite{mazzonetto2020drift}. However, a key challenge in our setting is controlling the probability of crossing a threshold between two consecutive observations. This is carefully addressed in the proof of Proposition~\ref{Hitting_time}, 
which we consider one of the main contributions of this paper, as it paves the way for extending our results to more general threshold diffusions. 
Moreover, to the best of our knowledge, Theorem~\ref{Thm_sigma}, which concerns diffusion coefficient estimation, is the first result of its kind in the context of threshold diffusions observed in high frequency and long time. However, in practice, high frequency data (fixed time horizon) is the usual setting for volatility estimation, see, e.g.,~\cite{ait2014high}, but obtaining similar results for T-CKLS would again require the knowledge of the law of the process. 
Theorem~\ref{Theorem_discret_cons_tcl} is the main result of this paper. Up to our knowledge, it is the first result treating parameter estimation for standard CKLS process from discrete observation. The same holds for T-CIR process and T-CKLS process. 
Moreover, we recover existing results on drift MLE for standard CIR and T-OU process providing a unified setting.  \\

In Section \ref{sec:numerics}, we evaluate the accuracy of our estimators and their asymptotic normality properties on a simulated dataset. Then, we apply these estimators to the ten-year US Treasury rate using the Federal Reserve Bank's H15 dataset. In Section~\ref{Interest rates analysis}, we use a discretized version of the test introduced in~\cite{su2017testing} to assess the presence of one or more thresholds in the drift term of the diffusion model. However, further theoretical development would be needed to rigorously establish the estimation of multiple thresholds, which is beyond the scope of this paper.

\paragraph{Outline.} In Section~\ref{Preliminaries}, first, we introduce the model, then we introduce the quasi-likelihood and likelihood functions associated to the T-CKLS process. We also comment on the assumptions of this article providing parameter sets for which all results of this paper hold. 
In Section~\ref{sec:continuous_th}, we deal with estimation of T-CKLS from continuous observations. We provide some statistical properties related to the drift estimator. 
The main results of this article are provided in Section~\ref{discret}, which deals with drift and volatility estimation from discrete observations. We study the asymptotic behavior in high frequency and long time for both estimators. 
Numerical experiments are provided in Section~\ref{sec:numerics}, where the estimators are implemented and tested on simulated data and US interest rates data. 
The proofs of our results are collected in Section~\ref{Proof}. Useful properties of the model are available in Section~\ref{appendix}.

Throughout the paper, we use the notion of stable convergence,  denoted $\xrightarrow[]{stably}$. 
Further details on this type of convergence can be found in \cite{jacod2011discretization} and \cite{jacod2013limit}.

\section{The framework: model and assumptions} \label{Preliminaries}

In the entire paper, $(\Omega, \mathcal{F}, (\mathcal{F}_t)_{t \geq 0}, \mathbb{P})$ denotes a filtered probability space satisfying the usual conditions and $B$ an $(\mathcal{F}_t)$-standard Brownian motion. 

\subsection{The model} \label{sec:toy}
The T-CKLS process solves the following one-dimensional SDE:
\begin{equation}
X_t=X_0+\int_0^t(a(X_s)-b(X_s)X_s)\vd s+\int_0^t\sigma(X_s)|X_s|^{\gamma(X_s)}\vd B_s, \quad t\geq 0,
\label{Multi_CKLS}
\end{equation}
with $X_0>0$, $X_0$ independent of $(B_t)_{t\geq 0}$, piecewise constant coefficients $a$, $b$, $\sigma$ and $\gamma$ possibly discontinuous at levels $0=r_0<r_1<\ldots<r_d<r_{d+1}=+\infty$, $d\in\mathbb{N}$. 
We focus on the case $\gamma(\mathbb{R}) \subseteq [0,1]$ and $\gamma(0) \in [{1}/{2},1] \cup \{0\}$. More precisely, let $I_j:=[r_j,r_{j+1})$, for $j\in\{0,\ldots,d\}$\footnote{One could also take  $(I_j)_{j= 0,\ldots, d}$ disjoint intervals such that $\cup_{j=0}^d I_j = [0,+\infty)$ and whose closure is $\overline I_j= [r_j,r_{j+1}]\cap \mathbb{R}$}, unless $\gamma(0)=0$, in which case $I_0=(-\infty,r_1)$. The drift coefficients are given by
\begin{equation}
a(x)=\sum_{j=0}^d a_j\mathds{1}_{I_j}(x) \in \mathbb{R} \quad \text{and} \quad b(x)=\sum_{j=0}^d b_j\mathds{1}_{I_j}(x)\in \mathbb{R},
\end{equation}
and similarly, the volatility coefficients are given by 
\begin{equation}
\sigma(x)=\sum_{j=0}^d \sigma_j \mathds{1}_{I_j}(x)>0 \quad \text{and} \quad \gamma(x)=\sum_{j=0}^d \gamma_j\mathds{1}_{I_j}(x) \geq 0.
\end{equation}
When $\gamma_0\in [1/2,1)$, we also assume that $a_0>0$.
When $\gamma_0=1$, we allow for $a_0 \geq 0$.
When $\gamma_0=0$ we keep the assumptions $r_1>0$ and $X_0>0$, although, in this specific case, one could allow for $X_0 \in \mathbb{R}$.

For existence of a pathwise unique strong solution to~\eqref{Multi_CKLS} (which is strong Markov) under the assumption that $\gamma_0\in[{1}/{2},1] \cup \{0\}$, we refer for instance to~\cite{engel_schm_89_3,le2006one}.
Moreover, it can be shown that T-CKLS is a Markov process (see, e.g.,~\cite{engel_schm_89_3}).
When $\gamma_0 \in [1/2,1]$ the process is always non-negative (see Lemma~\ref{lem:state_space} in Section~\ref{sec:property_stationary} for a precise statement). \\

We recall that in this article, we suppose that the thresholds $r=(r_j)_{j=0}^{d}$ and the diffusion exponents $\gamma=(\gamma_j)_{j=0}^d$ are known, and we estimate drift and diffusion parameter vectors $a,b,\sigma$ for continuous time observations and discrete (not necessarily equally spaced) high frequency observations and infinite horizon. We denote $\theta_{\star}:=(a,b)=(a_j,b_j)_{j=0}^d$ and $\sigma_\star$ the parameters to be estimated. 

The parameters for which our results hold are further restricted according to the following requirements: the likelihood (resp.~quasi-likelihood) function should be well defined (leading to the parameter space  $\Theta^{(\Lk)}$, resp.~$\Theta^{(\qLk)}$), the process should be ergodic, its stationary distribution should have finite moments of precise orders. 

In the remainder of this section we discuss these restrictions to derive the parameter space.
For the reader convenience, we now provide subsets of the state space for which all results of the paper hold: in addition to 
$\gamma(\mathbb{R}) \subseteq [0,1]$ and $ \sigma(\mathbb{R})\subseteq (0,+\infty)$, let  
$\gamma\equiv 0$ or let
\begin{equation} \label{eq:toy:noCIR}
    \gamma_0\in (1/2,1], \quad a_0>0, \quad b_d >0, 
\end{equation}
or let 
\begin{equation} \label{eq:toy:CIR}
    \gamma_0=1/2, \quad a_0 > \sigma_0^2, \quad b_d >0.
\end{equation} 
The latter case includes also standard CIR process treated in~\cite[Proposition 5]{ben2013asymptotic}, 
and 
CKLS treated in~\cite{mishura2022parameter} is included in \eqref{eq:toy:noCIR}. 
Another interesting (already treated in the literature) special case is $\gamma\equiv 0$ which corresponds either to T-OU process or T-drifted BM. 

\subsection{Likelihood and quasi-likelihood functions}

In the next sections, we assume $\sigma$ to be unknown, and we propose an estimator. Yet, in the following lines, the reader should think as if $\sigma$ is known (or replaced by an estimator). 

Assume that we have access to an observation of an entire trajectory on the time interval $[0,T]$ of the T-CKLS. 
We consider two different contrast functions: likelihood and a quasi-likelihood. 
The likelihood function $\theta \mapsto \Lk_T(\theta; \sigma, \gamma)$ is related to the Girsanov weight:
\begin{align}
\Lk_T(\theta; \sigma, \gamma)&=\exp\left(\int_0^T \frac{ a(X_s)-b(X_s)X_s}{\sigma(X_s)^2(X_s)^{2\gamma(X_s)}}\vd X_s -\frac{1}{2}\int_0^T\frac{\left(a(X_s)-b(X_s)X_s\right)^2}{\sigma(X_s)^2(X_s)^{2\gamma(X_s)}}\vd s\right).
    \label{Likelihood} 
\end{align}
Note that it is well defined if the integrals above are well defined (see Proposition \ref{Proposition_integral} for details). 
Hence, we restrict the parameter space $\Theta^{(\Lk)}$ to the coefficients for which the integrals above are well defined: if 
$\gamma_0={1}/{2}$, we have to restrict to $a_0\geq {\sigma_0^2}/{2}$.

We consider the quasi-likelihood function $\theta \mapsto \qLk_T(\theta):= \ln \Lk_T(\theta; 1, 0)$ \cite{su2015quasi}. Hence, 
\begin{equation}\label{Quasi_likelihood}
    \qLk_T(\theta)=\int_0^T \left( a(X_s)-b(X_s)X_s \right) \vd X_s -\frac{1}{2}\int_0^T\left(a(X_s)-b(X_s)X_s\right)^2 \vd s.
\end{equation}
The advantage of this contrast function is the fact that it does not depend on the diffusion's coefficients $\gamma,\sigma$. We denote $\Theta^{(\qLk)}$ the set of parameters such that the quasi-likelihood is well defined. Note that $\qLk_T$ is always well defined because the process we consider has continuous trajectories. So, it covers a wider range of parameters than the likelihood function.

\subsection{About our assumptions on ergodicity and moments of the stationary distribution} \label{sec:hypothesis:general}

We are interested only in the case in which the process is \emph{ergodic} or \emph{stationary}.
With \emph{ergodic}, we mean \emph{positive-recurrent}. When the process is ergodic, there exists a stationary distribution (invariant distribution for the transition semigroup), denoted by $\mu$. 
We say that the process is \emph{stationary} if it is ergodic, $X_0$ is independent of the driving Brownian motion, and $X_0$ is distributed according to the stationary distribution $\mu$.

The precise parameter restrictions ensuring ergodicity are provided in Table~\ref{table_ergo} in Section~\ref{sec:property_stationary}. 
An expression for $\mu$ for the T-CKLS is given in Section~\ref{sec:property_stationary} and its moments are studied in Proposition~\ref{Moment_sup}. By moment of order $m$, we mean $\int |x|^m \mu(\vd x) $. 

For clarity, the reader may keep the following cases in mind: 
The T-CKLS is ergodic if, for instance, the drift coefficients of the first interval $I_0$ and the last interval $I_d$ satisfy  
\[ 
    (a_0,b_0)\in (0,+\infty) \times \mathbb{R} \quad \text{ and } \quad (a_d,b_d)\in \mathbb{R} \times (0,+\infty)
\]
when $\gamma_0\neq 0$; and $a_0,a_d\in \mathbb{R}$ and $b_0, b_d>0$ when $\gamma_0=0$. Moreover, with the latter choices $(b_d\neq 0)$, the stationary distribution $\mu$ admits non-negative moments of all order  and, if $\gamma_0\not\in \{0,1/2\}$, negative moments of all order too. 
When $\gamma_0=0$ one can only have finite moments of negative order $m\in (-1,0)$. 
When $\gamma_0=1/2$, 
an additional parameter restriction on $a_0, \sigma_0$ is needed to ensure the existence of the moment of order $m<0$: $a_0> -m \sigma_0^2/2$. 

Our main results hold under some assumptions on the moments of the stationary distribution. 
In particular,  
\begin{equation}
     \mathbf{H}_{\qLk} = \calH_{\qLk} \colon 
    \mu \text{ admits finite } (2+2\gamma_d)\text{-th moment},
\end{equation}
which is  weaker than asking for the $4$-th order moment required in~\cite{su2015quasi}. 
The other two assumptions $\mathbf{H}_{\Lk}$ and $\calH_{\Lk}$ coincide with the previous if $\gamma\equiv 0$, otherwise they may involve negative moments:  

\begin{itemize}
    \item $\mathbf{H}_{\Lk}$: 
    $\mu$ admits finite $(-2\gamma_0)$-th and $(2-2\gamma_d)$-th moment, 
    \item $\calH_{\Lk}$: $\mu$ admits finite $(-q_{\Lk})$-th moment and $\max(2(1+\gamma_d), p_{\Lk})$-th moment with $p_{\Lk} > 1$ and $q_{\Lk} > 2$ such that $\frac{1}{p_{\Lk}} + \frac{2}{q_{\Lk}} = 1$.     
\end{itemize}
$\calH_{\Lk}$ is the most restrictive assumption. 
Since $p_{\Lk}=q_{\Lk}=3$ is a possible choice, we see that $\calH_{\Lk}$ is weaker than : $\mu$ admits moments of order $4$ and $-3$. 
Note that the case $\gamma_0=0$, $\gamma\not\equiv 0$ is not included in $\calH_{\Lk}$ because the process does not admit moments of any order smaller than than or equal than $-1$. 

We have seen that all assumptions $\mathbf{H}_{\Lk}, \mathbf{H}_{\qLk}=\calH_{\qLk}, \calH_{\Lk}$ are satisfied with the choices in~\eqref{eq:toy:noCIR} and \eqref{eq:toy:CIR}. 
As far as it concerns the other cases, the negative moments 
$(-q_{\Lk})$ and $(-2\gamma_0)$  impose restrictions on $a_0,\sigma_0$ only if $\gamma_0=1/2$. Therefore, the interplay between positive and negative moments $p_\Lk, q_\Lk$ is of interest in the case, $\gamma_0=1/2$. When $b_d = 0$ the condition on positive moments excludes some cases. The positive moment $2(1+\gamma_d)$ is not finite if $\gamma_d\in (1/2,1]$ and, if $\gamma_d=1/2$, it is finite only if $a_d<-3\sigma_d^2/2$. 
The positive moment $2(1-\gamma_d)$ is not finite if $\gamma_d\in [5/6,1)$ and imposes restrictions on the parameters if $\gamma_d=1/2$: $a_d<-\sigma_d^2/2$.

\section{Estimation from continuous time observations}\label{sec:continuous_th}

Let $T\in(0,\infty)$, and assume that we have at our disposal an entire trajectory of $(X_t)_{t\in [0,T]}$ 
of the 
solution to the SDE \eqref{Multi_CKLS}.

First, we provide estimators which maximize likelihood and quasi-likelihood. 
Next, we study the asymptotic behavior of the estimators in long time under the assumption that the process is ergodic (see Section~\ref{Preliminaries}).

\subsection{Estimators expressions}

The drift parameters estimators, MLE and QMLE, are defined as the maximizing argument of the log-likelihood \eqref{Likelihood} and the quasi-likelihood \eqref{Quasi_likelihood}: \begin{equation}
    \theta_T^{(\Lk)}:=\underset{\theta\in\Theta^{(\Lk)}}{\text{argmax }}\ln \Lk_T(\theta; \sigma, \gamma) 
    \quad \text{and} \quad 
    \theta_T^{(\qLk)}:=\underset{\theta\in\Theta^{(\qLk)}}{\text{argmax }}\qLk_T(\theta) .
\end{equation} 
We look for expressions for MLE and QMLE in terms of the following quantities:
\begin{equation}\label{Quantity_likelihood_continue_multi}
    \calQ_T^{j,m}:=\int_0^T X_s^m\mathds{1}_{I_j}(X_s)\vd s 
    \quad \text{and} \quad 
    \calM_T^{j,m}:=\int_0^T X_s^m\mathds{1}_{I_j}(X_s)\vd X_s
\end{equation}
for $j\in\{0,\ldots,d\}$ and $m\in\{-2\gamma_j,1-2\gamma_j,2-2\gamma_j,2\gamma_j\}\cup\{0,1,2\}$.
It is convenient to express log-likelihood and quasi-likelihood as follows.
The log-likelihood satisfies 
\begin{equation}
\ln \Lk_T(\theta; \sigma, \gamma):=\sum_{j=0}^d\frac{1}{\sigma_j^2} \left(a_j\calM_T^{j,-2\gamma_j}-b_j\calM_T^{j,1-2\gamma_j}-\frac{a_j^2}{2}\calQ_T^{j,-2\gamma_j}-a_j b_j\calQ_T^{j,1-2\gamma_j}-\frac{b_j^2}{2}\calQ_T^{j,2-2\gamma_j}\right) ,
\label{rewrite_L_cont}
\end{equation}
and the quasi-likelihood  
$\qLk_T(\theta) := \ln \Lk_T(\theta; 1, 0)$, which rewrites as
\begin{equation} \label{rewrite_Q_cont}
\qLk_T(\theta)=\sum_{j=0}^d a_j\calM_T^{j,0}-b_j\calM_T^{j,1}-\frac{a_j^2}{2}\calQ_T^{j,0}-a_j b_j\calQ_T^{j,1}-\frac{b_j^2}{2}\calQ_T^{j,2}.
\end{equation}

The following proposition provides explicit expression of MLE and QMLE, in terms of the quantities in~\eqref{Quantity_likelihood_continue_multi}.
\begin{prop}
\label{Proposition_estimateur_continu}
Let $T\in(0,\infty)$, the maximum of the likelihood $\Lk_T(\theta;\sigma,\gamma)$ is achieved at $\theta_T^{(\Lk)}:=(a_T^{j,\gamma},b_T^{j,\gamma})_{j=0}^d$ with
\begin{equation} \label{estimateurs_likelihood_continue_multi}
    (a_T^{j,\gamma},b_T^{j,\gamma}) =\left(\frac{\calM_T^{j,-2\gamma_j}\calQ_T^{j,2-2\gamma_j}-\calQ_T^{j,1-2\gamma_j}\calM_T^{j,1-2\gamma_j}}{\calQ_T^{j,-2\gamma_j}\calQ_T^{j,2-2\gamma_j}-(\calQ_T^{j,1-2\gamma_j})^2},\frac{\calM_T^{j,-2\gamma_j}\calQ_T^{j,1-2\gamma_j}-\calQ_T^{j,-2\gamma_j}\calM_T^{j,1-2\gamma_j}}{\calQ_T^{j,-2\gamma_j}\calQ_T^{j,2-2\gamma_j}-(\calQ_T^{j,1-2\gamma_j})^2}\right).
\end{equation}
The maximum of the quasi-likelihood $\qLk_T(\theta)$ is achieved at 
$\theta_T^{(\qLk)}:=(a_T^{j,0},b_T^{j,0})_{j=0}^d$, that is
\begin{equation} \label{estimateurs_quasi_continu_multi}
    (a_T^{j,0}, b_T^{j,0})=\left(\frac{\calM_T^{j,0}\calQ_T^{j,2}-\calQ_T^{j,1}\calM_T^{j,1}}{\calQ_T^{j,0}\calQ_T^{j,2}-(\calQ_T^{j,1})^2},\frac{\calM_T^{j,0}\calQ_T^{j,1}-\calQ_T^{j,0}\calM_T^{j,1}}{\calQ_T^{j,0}\calQ_T^{j,2}-(\calQ_T^{j,1})^2}\right).
\end{equation}
\end{prop}
\begin{proof}
We sketch the proof for MLE. The same works for QMLE.
One shows that~\eqref{estimateurs_likelihood_continue_multi} is the unique singular point of the gradient (vector of the derivatives with respect to $a_j$ and $b_j$ for all $j\in\{0,\ldots,d\}$) of~\eqref{rewrite_L_cont} and the Hessian is negative definite.
\end{proof}

\begin{remark} \label{rem:TOU}
If $\gamma\equiv 0$, then the diffusion coefficient is piecewise constant, so T-CKLS is a threshold OU (T-OU) and
QMLE and MLE coincide, as noticed in~\cite{mazzonetto2020drift}. 
\end{remark}

\begin{remark}
For every $j=0,\ldots,d$,
$(a_T^{j,\gamma},b_T^{j,\gamma})$ only depend on the observations of the trajectory $t\mapsto X_t$ which belong to $I_j$. Of course, the same holds for $(a_T^{j,0},b_T^{j,0})$.
\end{remark}

MLE and QMLE do not depend explicitly on $\sigma_\star$, but only on the quantities in~\eqref{Quantity_likelihood_continue_multi}. 
As usual, $\sigma_\star$ is related to the quadratic variation. 
Here, for every $j$, $\sigma_j^2$ is related to the quadratic variation of $\calM_T^{j,0}$, which is $\sigma_j^2 Q^{j,2\gamma_j}_T$.
By finding an a.s.~equivalent expression for $\langle \calM^{j,0} \rangle_T$ in terms of $\calM_T^{\cdot,0}, \calM_T^{j,1}, X_T, X_0$, we show in the next result that $\sigma_\star$ can be determined from a continuous trajectory. From this expression we derive the estimator from discrete observations defined in Theorem~\ref{Thm_sigma}.

\begin{prop}
\label{Estimateur_sigma_cont}
Let $T\in(0,\infty)$ and $j\in\{0,\ldots,d\}$. 
Then 
\begin{equation}
    \sigma_j=\sqrt{\frac{\frQ_T^{j} }{\calQ_T^{j,2\gamma_j}}} \quad \text{a.s. on the event } \{\calQ_T^{j,0}>0\}
    \label{eq: estim_vol_j}
\end{equation}
where 
\begin{equation}\label{discret_bracket2}
    \begin{split}
    \frQ_T^{0} & := \left(f_{0}(X_T)\right)^2- \left(f_{0}(X_0)\right)^2 + 2 \left(r_{1} \calM_T^{0,0}-\calM_T^{0,1}\right) - 2 r_1 \frak{f}_0,
    \\
    \frQ_T^{d} & := \left(f_{d}(X_T)\right)^2- \left(f_{d}(X_0)\right)^2 + 2 \left(r_{d} \calM_T^{d,0}-\calM_T^{d,1}\right),
    \end{split}
\end{equation}
with
$f_0(x)=x \mathds{1}_{I_0}(x) + r_1 \mathds{1}_{[r_1,+\infty)}(x) $, $f_d(x)=\mathds{1}_{I_d}(x) (x-r_d)$,
and $\frak{f}_0 = \min(X_T,r_{1})-\min(X_0,r_{1})$,
and for $j\in \{1,\ldots,d-1\}$:
\begin{equation}
    \frQ_T^{j}:= \left(f_{j}(X_T)\right)^2 - \left(f_{j}(X_0)\right)^2
    - 2\calM_T^{j,1}
    + 2 r_j\calM_T^{j,0}  + 2 (r_{j+1}-r_j) \big(\frak{f}_{j+1}  - \sum_{i=j+1}^d\calM_T^{i,0} \big)
    \label{discret_bracket1}
\end{equation}
with $f_{j}(x)= \mathds{1}_{I_j}(x) (x-r_j) + (r_{j+1}-r_j) \mathds{1}_{[r_{j+1},+\infty)}(x)$ and $\frak{f}_j=\max(X_T,r_{j})-\max(X_0,r_{j})$.
\end{prop}
\begin{proof}
For simplicity, we focus on the case $j\in\{1,\ldots,d-1\}$. The cases $j\in\{0, d\}$ follow similarly, as the proof requires only minor adjustments. On the event $\{\calQ_T^{j,0}>0\}$, the trajectories spend a positive amount of time in $I_j$ and so $\calQ_T^{j,2\gamma_j}$ does not vanish.  
Applying It\^o-Tanaka formula (see \cite[Chapter VI, exercice 1.25]{revuz2013continuous}) ensures that 
\begin{equation}
    \vd  f_{j}(X_s)
    = \mathds{1}_{I_j}(X_s) \vd X_s + 2^{-1} \vd \left( L^{r_j}_s(X)-L^{r_{j+1}}_s(X)\right).
\end{equation}
The quadratic variation of $f_{j}(X)$ satisfies a.s.~the equality:
\[
    \langle f_{j}(X) \rangle_T = \sigma_j^2 \int_0^T (X_s)^{2\gamma_j}\mathds{1}_{I_j}(X_s)\vd s =  \sigma_j^2 \calQ_T^{j,2\gamma_j}.
\]
Since for every semi-martingale $Y$, It\^o formula ensures that a.s.~$\vd Y^2 = 2 Y \vd Y + \vd\langle Y \rangle$, again It\^o formula applied to the semi-martingale $\left(f_{j}(X_T)\right)^2$, yields that a.s.~ 
\begin{equation}
    \vd \langle f_{j}(X) \rangle_s 
    = \vd \left(f_{j}(X_s)\right)^2 + 2 r_j \vd \calM_s^{j,0} - 2 \vd \calM_s^{j,1} + (r_{j+1} - r_{j}) \vd L^{r_{j+1}}_s(X).
\end{equation}
In order to check that 
$ \langle f_{j}(X) \rangle_T$
is a.s.~equal to~$\frQ_T^j$, we rewrite the local times in terms of $\calM^{j,0},\calM^{j,1}$:  
It\^o-Tanaka formula applied to $\max{(X_T,r_j)}$ yields 
$
    L^{r_j}_T(X) 
    = 2 \frak{f}_j - 2 \sum_{i=j}^d \calM_T^{i,0}
$
a.s.. The proof is thus completed.
\end{proof}

\subsection{Asymptotic properties: long time}
In this section, we explore the statistical properties as $T\to\infty$ of the MLE and QMLE from continuous time observations of a trajectory of the T-CKLS process. 

We assume that the process is ergodic, $\mu$ is the stationary distribution, and we introduce the following hypotheses:
\begin{itemize}
   \item $\mathbf{H}_{\Lk}$: $\mu$ admits finite $(-2\gamma_0)$-th and $(2-2\gamma_d)$-th moment, 
    \item $\mathbf{H}_{\qLk}$: $\mu$ admits finite $(2+2\gamma_d)$-th moment.
\end{itemize}
We have discussed these assumptions in Section~\ref{sec:hypothesis:general}, and we have provided a set of parameters for which the above assumptions are satisfied in Section~\ref{sec:toy}. 

The asymptotic behavior of the MLE and QMLE are provided in the following theorem, which states that MLE and QMLE are strongly consistent and asymptotically normal estimators of the drift parameter $\theta_{\star}$.

\begin{theorem}
\label{Theorem_cont}
For \(\ell \in \{\Lk, \qLk\}\), under Hypothesis \(\mathbf{H}_\ell\),  $\theta_T^{(\ell)}$ (MLE $\theta_T^{(\Lk)}$ or QMLE $ \theta_T^{(\qLk)}$) is a strongly consistent estimators of \(\theta_{\star}\), \textit{i.e.}
\begin{equation}
    \theta_T^{(\ell)} \xrightarrow[T \to +\infty]{a.s.} \theta_{\star}.
    \label{Consitance_Cont}
\end{equation}
Additionally, the following convergence holds:
\begin{equation}\label{TCL_Cont}
     \sqrt{T}\left(\theta_T^{(\ell)} - \theta_{\star}\right) \xrightarrow[T \to +\infty]{\text{stably}} \mathbf{N}^{(\ell)},
\end{equation}
where \(\mathbf{N}^{(\ell)} = (N_{j,a}^{(\ell)}, N_{j,b}^{(\ell)})_{j=0}^d\) are \(d+1\) independent two-dimensional centered Gaussian random variables, independent of \(X\), with covariance matrices respectively given by \(\sigma_j^2 \Gamma_j^{(\ell)}\):
\[
    \Gamma_j^{(\Lk)} := \Gamma_j^{(\Lk,\gamma_j)} := 
    \begin{pmatrix}
        \calQ_\infty^{j,-2\gamma_j} & -\calQ_\infty^{j,1-2\gamma_j} \\
        -\calQ_\infty^{j,1-2\gamma_j} & \calQ_\infty^{j,2-2\gamma_j}
    \end{pmatrix}^{\!\!-1}
    \quad \text{and} \quad
    \Gamma_j^{(\qLk)} = 
    \left(\Gamma_j^{(\Lk,0)}\right)^{-1}
    \Gamma_j^{(\Lk,-\gamma_j)}   
    \left(\Gamma_j^{(\Lk,0)}\right)^{-1},
\]
where \(\calQ_\infty^{j,.}\) are real constants defined in the next Lemma \ref{Quantité_regime_multi}. 

Furthermore, the LAN property holds for the \(\ell\)-function with a convergence rate of ${1}/{\sqrt{T}}$, and the asymptotic Fisher information is given by
\begin{align}
    & \Gamma^{(\Lk)}_{\text{LAN}} = \text{diag}\left( \frac{1}{\sigma_0^2}  (\Gamma_0^{(\Lk)})^{-1} , \frac{1}{\sigma_1^2} (\Gamma_1^{(\Lk)})^{-1} , \dots, \frac{1}{\sigma_d^2}  (\Gamma_d^{(\Lk)})^{-1}  \right)\quad \\ 
    & \text{and}\quad 
    \Gamma^{(\qLk)}_{\text{LAN}} = \text{diag}\left( \sigma_0^2  (\Gamma_0^{(\Lk,-\gamma_j)}   )^{-1} , \sigma_1^2  (\Gamma_1^{(\Lk,-\gamma_j)}   )^{-1} , \dots, \sigma_d^2  (\Gamma_d^{(\Lk,-\gamma_j)}   )^{-1}  \right).
\end{align}
\end{theorem}

\begin{proof}
The proof is similar to those in~\cite{lejay2020maximum, mazzonetto2020drift}, so we summarize the steps while emphasizing the specifics of the case considered in this article. Note that MLE and QMLE can be rewritten as follows: 
\[
    \theta_j^{(\Lk)}
    \overset{a.s.}{=}
    \theta_j
    +\sigma_j     
    \begin{pmatrix}
    \mathfrak{M}_T^{j,-\gamma_j} 
    ,
    -\mathfrak{M}_T^{j,1-\gamma_j}
    \end{pmatrix}
    \Gamma^{(\Lk)}
    \quad \text{and} \quad
    \theta_j^{(\qLk)}
    \overset{a.s.}{=}
    \theta_j
    +\sigma_j     
    \begin{pmatrix}
    \mathfrak{M}_T^{j,\gamma_j} 
    ,
    -\mathfrak{M}_T^{j,1+\gamma_j}
    \end{pmatrix}
    \Gamma^{(\Lk,0)}
\]
where 
$\mathfrak{M}_T^{j,k}=\int_0^T (X_s)^{k}\mathds{1}_{I_j}(X_s)\vd B_s$ 
for $k\in\{-\gamma_j,1-\gamma_j,\gamma_j,1+\gamma_j\}$
are martingales. Indeed, note that $\calQ_T^{j,2k}$ is the quadratic variation of $\calM_T^{j,k-\gamma_j}$ and of $\mathfrak{M}_T^{j,k}$ up to a multiplicative factor.  
The consistency of the estimator (MLE and QMLE) follows directly from~\cite[Theorem 1]{lepingle2006comportement} 
and the ergodicity of the process which implies Lemma~\ref{Quantité_regime_multi}. 
The asymptotic normality property follows from the martingale convergence theorem~\cite[Theorem 2.2]{crimaldi2005convergence}. 
Note that Hypotheses $\mathbf{H}_\Lk$ and $\mathbf{H}_{\qLk}$ are necessary for the application of~\cite[Theorem 1]{lepingle2006comportement}.

\textit{Proof of the LAN property.} 
Let $\Delta \theta_\star:=(\Delta a_j,\Delta b_j)_{j=0}^d$. Then, 
\begin{align}
    &\ln\frac{\Lk_T\left(\theta_\star+\frac{1}{\sqrt{T}}\Delta\theta_\star\right)}{\Lk_T(\theta_\star)}\\
    & =\sum_{j=0}^d\left(\frac{1}{\sigma_j\sqrt{T}}\begin{pmatrix}
    \Delta a_j, \Delta b_j
    \end{pmatrix}A_T^{j,(\Lk)} -\frac{1}{2\sigma_j^2T}\begin{pmatrix}
    \Delta a_j,
    \Delta b_j
    \end{pmatrix}\langle A^{j,(\Lk)}\rangle_T \begin{pmatrix}
    \Delta a_j,
    \Delta b_j
    \end{pmatrix}'\right),
\end{align}
and
\begin{align}
    &\qLk_T\left(\theta_\star+\frac{1}{\sqrt{T}}\Delta\theta_\star\right)-\qLk_T(\theta_\star)\\
    &=\sum_{j=0}^d\left(\frac{\sigma_j}{\sqrt{T}}\begin{pmatrix}
    \Delta a_j, \Delta b_j
    \end{pmatrix}A_T^{j,(\qLk)} 
    -\frac{\sigma_j^2}{2T}\begin{pmatrix}
    \Delta a_j,
    \Delta b_j
    \end{pmatrix}\langle A^{j,(\qLk)}\rangle_T\begin{pmatrix}
    \Delta a_j,
    \Delta b_j
    \end{pmatrix}'\right),
\end{align}
with $A_T^{j,(\Lk)}=\left(\mathfrak{M}_T^{-\gamma_j},-\mathfrak{M}_T^{1-\gamma_j} \right)'$ and $A_T^{j,(\qLk)}=\left(\mathfrak{M}_T^{\gamma_j},-\mathfrak{M}_T^{1+\gamma_j} \right)'$. Then, for $\ell\in\{\Lk,\qLk\}$, by applying \cite[Theorem 2.2]{crimaldi2005convergence} and Lemma \ref{Quantité_regime_multi}, we have $\langle A^{j,(\ell)}\rangle_T/T\xrightarrow[T\to +\infty]{a.s.}(\Gamma_j^{(\ell)})^{-1}$ and $(A_t^{j,(\ell)})_{j=0}^d /{\sqrt{T}}\xrightarrow[T\to\infty]{law}\mathcal{N}(0,\Gamma^{(\ell)}_{\text{LAN}})$. This completes the proof.
\end{proof}

\begin{lemma}[Ergodic properties]
\label{Quantité_regime_multi}
   For $j\in\{0,\ldots, d\}$ and $m\in\mathbb{R}$, if the $m$-th moment of $\mu$ is finite on the set $I_j$, then $$\calQ_\infty^{j,m}\overset{a.s.}{:=}\lim_{T\to \infty} \frac{\calQ_T^{j,m}}{T}=\int_{I_j}x^m\mu(\vd x),$$
    are non-vanishing constants.    
\end{lemma}

\begin{remark} \label{coroll:LAN_continuous}
The LAN property, is a fundamental concept in the asymptotic theory of statistics. For instance, when it is satisfied, it can be combined with the Minimax theorem to establish a lower bound for the asymptotic variance of estimators. For more details, we refer to the book \cite{le2000asymptotics}.
\end{remark}

\begin{remark}[Estimator of the asymptotic variance] 
\label{Remark:Var_cons}
We can obtain a consistent estimator for the asymptotic variance using the following construction. Define:  
\[
    \Gamma_{j,T}^{(\Lk)} := \Gamma_{j,T}^{(\Lk,\gamma_j)} := 
    \frac{1}{T}\begin{pmatrix}
        \calQ_T^{j,-2\gamma_j} & -\calQ_T^{j,1-2\gamma_j} \\
        -\calQ_T^{j,1-2\gamma_j} & \calQ_T^{j,2-2\gamma_j}
    \end{pmatrix}^{\!\!-1}\quad \text{and} \quad \Gamma_{j,T}^{(\qLk)} := 
    \left(\Gamma_{j,T}^{(\Lk,0)}\right)^{-1}  
    \Gamma_{j,T}^{(\Lk,-\gamma_j)}   
    \left(\Gamma_{j,T}^{(\Lk,0)}\right)^{-1}.
\]    
Using these, we can construct $\Gamma^{(\ell)}_T$. By applying Lemma \ref{Quantité_regime_multi}, it follows that:  
\begin{equation}
    \Gamma^{(\ell)}_T \xrightarrow[T \to +\infty]{a.s.} \Gamma^{(\ell)}.
\end{equation}  
\end{remark}

\section{Estimation from discrete observations}\label{discret}

In this section, we assume we observe the process, solution to the SDE \eqref{Multi_CKLS}, on a discrete time grid $0=t_0<t_1<\ldots<t_{N-1}<t_N=T_N<\infty$, for $N\in\mathbb{N}$. 

First, we provide estimators, which maximize a discretized versions of likelihood  
and quasi-likelihood. 
Next, we study the asymptotic behavior of the estimators in high frequency and long time under the assumption that the process is stationary (see Section~\ref{Preliminaries}).

\subsection{Estimators expressions} 
 \label{sec:estimators:discrete}
There is no exploitable explicit expression for the finite dimensional distributions of the T-CKLS process.
Hence, instead of considering the likelihood function associated to the sample $(X_{t_k})_{k=0}^N$, we considered a discretization of the likelihood $\Lk_T$~\eqref{Likelihood} and quasi-likelihood $\qLk_T$~\eqref{Quasi_likelihood}. 
Once these discretizations introduced, we denote them respectively by~$\Lk_{T_N,N}$ and $\qLk_{T_N,N}$ and
we compute the estimators 
\begin{equation}
    \theta_{T_N,N}^{(\Lk)}=\underset{\theta\in\Theta^{(\Lk)}}{\text{argmax } } \ln \Lk_{T_N,N}(\theta; \sigma, \gamma) 
    \quad \text{and} \quad
    \theta_{T_N,N}^{(\qLk)}=\underset{\theta\in\Theta^{(\qLk)}}{\text{argmax }}\qLk_{T_N,N}(\theta).
\end{equation}

Let us denote by $\calQ_{T_N,N}^{j,m}$ and $\calM_{T_N,N}^{j,m}$ the discrete versions of $\calQ_{T_N}^{j,m}$ and $\calM_{T_N}^{j,m}$ 
in~\eqref{Quantity_likelihood_continue_multi}:
\begin{equation}
\label{Quantity_likelihood_discrete_multi}    
\calQ_{T_N,N}^{j,m}:=\sum_{i=0}^{N-1}X_{t_i}^m\mathds{1}_{I_j}(X_{t_i})(t_{i+1}-t_i) \quad \text{and} \quad \calM_{T_N,N}^{j,m}:=\sum_{i=0}^{N-1}  X_{t_i}^m\mathds{1}_{I_j}(X_{t_i})(X_{t_{i+1}}-X_{t_i})
\end{equation}
for $j\in\{0,\ldots,d\}$ and $m\in\{-2\gamma_j,1-2\gamma_j,2-2\gamma_j,2\gamma_j\}\cup\{-1,0,1,2\}$.
To obtain the discretized likelihood and quasi-likelihood, it would be natural to replace the above quantities in the continuous time observations likelihood and quasi-likelihood functions given in~\eqref{rewrite_L_cont}-\eqref{rewrite_Q_cont}.
The discretized quasi-likelihood is then
\begin{equation} \label{rewrite_Q_di}
\qLk_{T_N,N}(\theta):=\sum_{j=0}^d a_j\calM_{T_N,N}^{j,0}-b_j\calM_{T_N,N}^{j,1}-\frac{a_j^2}{2}\calQ_{T_N,N}^{j,0}-a_j b_j\calQ_{T_N,N}^{j,1}-\frac{b_j^2}{2}\calQ_{T_N,N}^{j,2},
\end{equation}
and we \emph{could} do similarly for the log-likelihood. We would get
\begin{align}
\label{rewrite_L_di_0}
& \ln \Lk_{T_N,N}(\theta; \sigma, \gamma)
\\ & =\sum_{j=0}^d\frac{1}{\sigma_j^2} \left(a_j\calM_{T_N,N}^{j,-2\gamma_j}-b_j\calM_{T_N,N}^{j,1-2\gamma_j}-\frac{a_j^2}{2}\calQ_{T_N,N}^{j,-2\gamma_j}-a_j b_j\calQ_{T_N,N}^{j,1-2\gamma_j}-\frac{b_j^2}{2}\calQ_{T_N,N}^{j,2-2\gamma_j}\right).
\end{align}
Actually, we do not choose the latter quantity as discretized-log-likelihood. Instead, we consider the following discretized log-likelihood:
\begin{equation}\label{rewrite_L_di}
    \ln \Lk_{T_N,N}(\theta; \sigma, \gamma)=\sum_{j=0}^d\frac{1}{\sigma_j^2} \left(a_j\frM_{T_N,N}^{j,-2\gamma_j}-b_j\frM_{T_N,N}^{j,1-2\gamma_j}-\frac{a_j^2}{2}\calQ_{T_N,N}^{j,-2\gamma_j}-a_j b_j\calQ_{T_N,N}^{j,1-2\gamma_j}-\frac{b_j^2}{2}\calQ_{T_N,N}^{j,2-2\gamma_j}\right),
\end{equation}
which involves $\frM_{T_N,N}^{j,m}$, 
a different discretization of $\calM^{j,m}_{T_N}$. 
This new discretization is based on an alternative expression of $\calM^{j,m}_{T_N}$ (given in Lemma~\ref{lemma_disc} in Section~\ref{Proof}) which depends on $\calM^{j,0}_{T_N}$ and $\calQ^{j,0}_{T_N}$ and $\calQ^{j,-1}_{T_N}$ and it is obtained by replacing the last quantities by their discrete analogues given above.  
Taking the new discretization, for which we provide the expression in some lines, yields the assumptions on the negative moments of the stationary distribution in $\calH_{\Lk}$, which would be stronger otherwise. 
Indeed, the asymptotic behavior of $\calM_{T_N,N}^{0,m}$ is related to the one of $\calQ_{T_N,N}^{0,2m}$ and to its expectation, which is finite if $\mu$ has finite $2m$-order moments. 
If, for instance, $\gamma_0 =1/2$ the new discretization involves only $\calQ_{T_N,N}^{0,k}$ with $k=-1,0,1$ and $\calM_{T_N,N}^{0,0}$, while old one involves also $\calM_{T_N,N}^{0,-1}$. 
The idea of modifying the likelihood expression to obtain a more convenient form is inspired by~\cite{ben2013asymptotic}.

Let us now introduce these new discretizations, for every $j\in\{0,\ldots,d\}$ we define
$\frM^{j,0}_{T_N,N} := \calM^{j,0}_{T_N,N}$ 
and for $m\in\{-2\gamma_j,1-2\gamma_j\}\setminus\{0\}$, then
\begin{equation}
    \frM^{0,m}_{T_N,N}=  f_{0,m+1}(X_0)-f_{0,m+1}(X_{T_N})-\frac{m}{2}\sigma_0^2\calQ_{T_N,N}^{0,m+2\gamma_0-1}   + r_1^m\left(\calM_{T_N,N}^{0,0} +\frak{f}_0 \right),
    \label{frak_M_2_0_disc}
\end{equation}
\begin{equation}
    \frM^{d,m}_{T_N,N}=  f_{d,m+1}(X_{T_N})-f_{d,m+1}(X_0)-\frac{m}{2}\sigma_d^2\calQ_{T_N,N}^{d,m+2\gamma_d-1}   + r_{d}^m\left(\calM_{T_N,N}^{d,0} - \frak{f}_d \right),
    \label{frak_M_2_d_disc}
\end{equation}
 where $\frak{f}_0 = \min(X_T,r_{1})-\min(X_0,r_{1})$,
$f_{0,m}(x)=\int_{x}^{r_{1}} y^{m-1}dy \, \mathds{1}_{I_0}(x)$,
$f_{d,m}(x)=\int_{r_{d}}^x y^{m-1}dy \, \mathds{1}_{I_d}(x)$,
and if $j\in\{1,\ldots, d-1\}$, then
\begin{align}
\frM^{j,m}_{T_N,N} &= f_{j,m+1}(X_{T_N}) - f_{j,m+1} (X_0)-\frac{m}{2}\sigma_j^2\calQ_{T_N,N}^{j,m+2\gamma_j-1}\nonumber \\
    & \quad  + r_j^m \calM_{T_N,N}^{j,0} + r_{j+1}^m \frak{f}_{j+1}- r_{j}^m \frak{f}_{j} - (r_{j+1}^m - r_j^m) \sum_{k=j+1}^{d} \calM_{T_N,N}^{k,0}
    \label{frak_M_1_discr}
\end{align}
where $f_{j,m}(x)=\int_{r_j}^{x \wedge r_{j+1}}y^{m-1}dy \, \mathds{1}_{(r_j,+\infty)}(x)$ and $\frak{f}_j=\max(X_T,r_{j})-\max(X_0,r_{j})$.

The following proposition establishes an explicit expression of the discretized MLE and discretized QMLE.
\begin{prop}
\label{prop:discrete_estimator}
Let $(T_N)_{N\in\mathbb{N}}$ be a sequence in $(0,\infty)$ and let $N\in\mathbb{N}$.
The maximum of the discretized likelihood is achieved at $\theta_{T_N,N}^{(\Lk)}=(a_{T_N,N}^{j,\gamma_j},b_{T_N,N}^{j,\gamma_j})_{j=0}^d$ with
\begin{equation}
   (a_{T_N,N}^{j,\gamma_j},b_{T_N,N}^{j,\gamma_j}):=\left(\tfrac{\frM_{T_N,N}^{j,-2\gamma_j}\calQ_{T_N,N}^{j,2-2\gamma_j}-\calQ_{T_N,N}^{j,1-2\gamma_j}\frM_{T_N,N}^{j,1-2\gamma_j}}{\calQ_{T_N,N}^{j,-2\gamma_j}\calQ_{T_N,N}^{j,2-2\gamma_j}-(\calQ_{T_N,N}^{j,1-2\gamma_j})^2},
   \tfrac{\frM_{T_N,N}^{j,-2\gamma_j}\calQ_{T_N,N}^{j,1-2\gamma_j}-\calQ_{T_N,N}^{j,-2\gamma_j}\frM_{T_N,N}^{j,1-2\gamma_j}}{\calQ_{T_N,N}^{j,-2\gamma_j}\calQ_{T_N,N}^{j,2-2\gamma_j}-(\calQ_{T_N,N}^{j,1-2\gamma_j})^2} \right).
    \label{Estim_discret_multi_MLE}
\end{equation}
The maximum of discretized quasi-likelihood is achieved at $\theta_{T_N,N}^{(\qLk)}=(a_{T_N,N}^{j,0},b_{T_N,N}^{j,0})_{j=0}^d$ with
\begin{equation}
(a_{T_N,N}^{j,0},b_{T_N,N}^{j,0})=\left(\tfrac{\calM_{T_N,N}^{j,0}\calQ_{T_N,N}^{j,2}-\calQ_{T_N,N}^{j,1}\calM_{T_N,N}^{j,1}}{\calQ_{T_N,N}^{j,0}\calQ_{T_N,N}^{j,2}-(\calQ_{T_N,N}^{j,1})^2},\tfrac{\calM_{T_N,N}^{j,0}\calQ_{T_N,N}^{j,1}-\calQ_{T_N,N}^{j,0}\calM_{T_N,N}^{j,1}}{\calQ_{T_N,N}^{j,0}\calQ_{T_N,N}^{j,2}-(\calQ_{T_N,N}^{j,1})^2} \right).
    \label{Estim_discret_multi_QMLE_1}
\end{equation}
\end{prop}
The proof is omitted because it is analogous to the one of Proposition~\ref{Proposition_estimateur_continu}.

Note that the estimators in Proposition~\ref{prop:discrete_estimator} are the discretisations of the estimators in Proposition~\ref{Proposition_estimateur_continu}. 

Observe also that the QMLE does not depend on the parameter vector $\sigma_\star=(\sigma_j)_{j=0}^d$, instead the MLE does because so do the expressions $\frM_{T_N,N}^{j,m}$.
Since we assume $\sigma_\star$ is not known, we replace it by an estimator in the MLE. 
For $j\in\{0,\ldots,d\}$,
\begin{equation}\label{estim_sigma_discret}
\sigma^j_{T_N,N}=\sqrt{\frac{   \frQ^j_{T_N,N} }{\calQ_{T_N,N}^{j,2\gamma_j}}},
\end{equation}
where $\frQ^j_{T_N,N}$ is obtained by discretizing the right-hand side of formulas~\eqref{discret_bracket2} and \eqref{discret_bracket1} in Proposition~\ref{Estimateur_sigma_cont}. 
Note that $\frQ^j_{T_N,N}$ depends on $\calM^{j,1}_{T_N,N}$ and $\calM^{j,0}_{T_N,N}$ defined in~\eqref{Quantity_likelihood_discrete_multi}. 
For instance, 
\begin{equation}
    \frQ^d_{T_N,N} := \mathds{1}_{I_d}(X_T) (X_T-r_d)^2- \mathds{1}_{I_d}(X_0)  (X_0-r_d)^2 +2\left(r_d\calM_{T_N,N}^{d,0}-\calM_{T_N,N}^{d,1}\right).
\end{equation}

\begin{remark}[QMLE and $\frM^{j,1}_{T_N,N}$]
    Replacing $\calM^{j,1}_{T_N,N}$ with $\frM^{j,1}_{T_N,N}$ in the QMLE estimator does not relax the moment assumptions, as 
    the estimator~\eqref{estim_sigma_discret} of $\sigma_\star$ depends on $\calM^{j,1}_{T_N,N}$.  
    Likewise, if $\sigma_\star$ is known, there is no improvement either; see Remark~\ref{rem:ass:QMLE} for more details.  
\end{remark}

\subsection{Asymptotic properties: high frequency - long time}\label{sec:asymp_discret}
In this section, we state the statistical properties of the discretized MLE and QMLE. 

Let $\Delta_N:=\max_{k=1,\ldots,N}(t_k-t_{k-1})$ denote the maximal lag between two consecutive observations.
We assume that the observation time window goes to infinity (long time) and the maximal lag between consecutive observations vanishes  (high frequency):
\begin{equation} \label{eq:Lt-Hf}
    \lim_{N\to+\infty}T_N=+\infty \quad \text{and} \quad \lim_{N\to+\infty}\Delta_N=0. 
\end{equation}
We assume that the process is stationary, $X_0 \sim \mu$, where $\mu$ is the stationary distribution. 
The results of this section require additional assumptions on the moments of $\mu$.  
We start with the assumptions on the moments. If $\gamma\equiv 0$ then $\calH_{\Lk}=\calH_{\qLk}=\mathbf{H}_{\Lk}=\mathbf{H}_{\qLk}$ : $\mu$ admits finite second moment. If $\gamma\not\equiv 0$: 
\begin{itemize}
    \item $\calH_{\Lk}$ : 
   $\mu$ admits finite $\max(2(1+\gamma_d), p_{\Lk})$-th moment and $(-q_{\Lk})$-th moment, with $p_{\Lk} > 1$ and $q_{\Lk} > 2$ such that $\frac{1}{p_{\Lk}} + \frac{2}{q_{\Lk}} = 1$.     
    \item $\calH_{\qLk} = \mathbf{H}_{\qLk}$:  
    $\mu$ admits finite $(2+2\gamma_d)$-th moment.
\end{itemize}

We have discussed these assumptions in Section~\ref{sec:hypothesis:general} and, in Section~\ref{sec:toy}, we have provided a set of parameters for which the assumptions are satisfied. 

\subsubsection{Results} 
We are now ready to provide our first convergence result in the discrete setting.
We consider the volatility estimator in~\eqref{estim_sigma_discret} 
and we prove consistency, and show that the speed of convergence is larger than $\sqrt{T_N}$.

 \begin{theorem}\label{Thm_sigma}
 Assume that~\eqref{eq:Lt-Hf} holds, that the T-CKLS $X$ is stationary and that Hypothesis $\calH_{\qLk}$ holds. 
 Then, the estimator $\sigma_{T_N,N}^2=((\sigma_{T_N,N}^j)^2)_{j=0}^d$ in~\eqref{estim_sigma_discret} 
 is a consistent estimators of  the diffusion coefficient vector 
 $\sigma^2_\star= ( (\sigma_j)^2 )_{j=0}^d$, \textit{i.e.}
\begin{equation}
\sigma_{T_N,N}^2 \xrightarrow[N\to +\infty]{\mathbb{P}}\sigma^2_\star.
\label{consistance_estimateur_sigma}
\end{equation}
Under the additional assumption that $\underset{N\to+\infty}{\lim}T_N \Delta_N=0$, it holds that 
\begin{equation}
    \sqrt{T_N}\left(\sigma_{T_N,N}^2-\sigma^2_\star \right)\xrightarrow[N\to +\infty]{\mathbb{P}}0_{\mathbb{R}^{d+1}}.
\end{equation}
\end{theorem}
For classical diffusions, it is well known (see, e.g.,~\cite{ait2014high}) that in high frequency for fixed time horizon the estimator constructed using quadratic variation, such as ours, will achieve a convergence rate of $1 / \sqrt{\Delta_N}$ when $\Delta_N \to 0$ towards $\sigma_\star^2$. The same holds in the context of TDs with $\Delta_N=1/N$, see~\cite{lejay2018statistical}. Although we are in a different setting (long time), we expect a similar convergence rate; nevertheless, it seems that deeper knowledge of the law is required. 

Since we assume that the diffusion coefficient vector $\sigma_\star$ is unknown, 
we replace it by estimator~\eqref{estim_sigma_discret} in the expression of $\frM_{T_N,N}^{j,\cdot}$.

The following theorem states the asymptotic properties in high frequency and long time observations of the discretized MLE and QMLE of $\theta_{\star}$.     

\begin{theorem} 
\label{Theorem_discret_cons_tcl}
Assume that~\eqref{eq:Lt-Hf} holds and that the T-CKLS $X$ is stationary.
For $\ell\in\{\Lk,\qLk\}$, under Hypothesis $\calH_\ell$, the estimator $\theta_{T_N,N}^{(\ell)}$ (MLE $\theta_{T_N,N}^{(\Lk)}$ or QMLE $\theta_{T_N,N}^{(\qLk)}$) is a consistent estimator of $\theta_{\star}$ \text{i.e.}
    \begin{equation}\label{Constistance_discret}
        \theta_{T_N,N}^{(\ell)} \xrightarrow[N\to +\infty]{\mathbb{P}}\theta_{\star}.
    \end{equation}
Assume that $\underset{N\to+\infty}{\lim}T_N \Delta_N=0$. 
Then under the same hypothesis, we have
\begin{equation}\label{TCL_Discrete}
    \sqrt{T_N}\left(\theta_{T_N,N}^{(\ell)}-\theta_{\star}\right)\xrightarrow[N\to +\infty]{\text{stably}}\mathbf{N}^{\left(\ell\right)},
\end{equation}
where $\mathbf{N}^{(.)}$ is defined in Theorem \ref{Theorem_cont}. Additionally, the analog of the LAN property for the $\ell$-function holds at the true drift parameters
with a convergence rate of  $1/{\sqrt{T_N}}$ and an asymptotic Fisher information given as in Theorem \ref{Theorem_cont}.
\end{theorem}

The proof of Theorem \ref{Thm_sigma} and Theorem \ref{Theorem_discret_cons_tcl} are postponed to Section~\ref{Proof}. Nevertheless we give the main ideas and tools in Section~\ref{ssec:discussions:proofkey}.

\begin{remark}
    The above theorem holds whether $\sigma_{\star}$ is known or not known and replaced by its estimator~\eqref{estim_sigma_discret}.  
\end{remark}

\begin{remark}[Joint convergence]
For $\ell\in\{\Lk,\qLk\}$, it follows from the previous results that the estimator $(\theta_{T_N,N}^{(\ell)}, \sigma_{T_N,N})$ is a consistent estimator of $(\theta_{\star},\sigma_\star)$. By the stable convergence properties, 
$\sqrt{T_N}(\theta_{T_N,N}^{(\ell)}, \sigma_{T_N,N})$ converges stably to the vector $(\mathbf{N}^{\left(\ell\right)}, 0_{\mathbb{R}^{d+1}} )$. 
\end{remark}

\begin{remark}[Asymptotic variance estimation]  
\label{Remark:Var_disc}
Similar to Remark \ref{Remark:Var_cons}, a consistent estimator of the asymptotic variance can be constructed using a consequence of the combination of Lemma~\ref{Quantité_regime_multi} and the next Lemma~\ref{Convergence_Tps_long}:  
$    
    \calQ_{T_N,N}^{j,m} /T_N\xrightarrow[N \to +\infty]{\mathbb{P}} \calQ_\infty^{j,m},
$   
for suitable $m$.   
\end{remark}

\subsection{Key elements of the proofs} \label{ssec:discussions:proofkey}

We now comment on the proofs of consistency and asymptotic normality for drift and diffusion coefficients, namely Theorems~\ref{Thm_sigma} and \ref{Theorem_discret_cons_tcl}.
They rely on the next Lemma~\ref{Convergence_Tps_long} and on the continuous time results Proposition~\ref{Estimateur_sigma_cont} and 
Theorem~\ref{Theorem_cont}. 
More precisely, one shows that 
the rescaled difference between the discrete time and continuous time estimators 
(e.g.,~$\sqrt{T_N} \left(\theta^{(\Lk)}_{T_N,N} - \theta^{(\Lk)}_{T_N}\right)$)
vanishes as $N\to\infty$. Similarly, for the LAN property, the key idea is to control the difference between the log-likelihood ratios in the discrete and continuous settings, ensuring that their rescaled difference vanishes asymptotically.  
Since all estimators depend on $\calM^{j,m}_{T_N,N}$ and $\calQ^{j,k}_{T_N,N}$ for suitable $k,m$, the proofs are based on the following result. 

\begin{lemma}  \label{Convergence_Tps_long}
Let $\lambda\in\{1,2\}$. Assume that $\underset{N\to+\infty}{\lim}T_N=+\infty$, $\underset{N\to+\infty}{\lim}T_N^{\lambda-1} \Delta_N=0$, and that the T-CKLS $X$ is stationary (see Section~\ref{Preliminaries}).
Then
\begin{equation}\label{Conv_i_m}
    \lim_{N\to \infty} {T_N^{-\sfrac{1}{\lambda}}} \Emu\left[|\calQ_{T_N}^{j,k}-\calQ_{T_N,N}^{j,k}|\right]=0 
    \quad \text{and} \quad 
    \lim_{N\to \infty} {T_N^{-\sfrac{1}{\lambda}}}\Emu\left[|\calM_{T_N}^{j,m}-\calM_{T_N,N}^{j,m}|\right]=0
\end{equation}
for all $j\in\{0,\ldots,d\}$ in each one of the following cases
\begin{enumerate}[(a)]
  	\item \label{item_1}
    (drift QMLE)
    for every $k\in\{0,1,2\}$ and $m\in\{0,1\}$,
    under hypothesis $\calH_{\qLk}$.
    \item \label{item_0} 
    (volatility estimation)
    for every $k=2\gamma_j$ and $m\in\{0,1\}$,
    under hypothesis $\calH_{\qLk}$. 
      \item \label{item_2} 
    (drift MLE, $\sigma_\star$ unknown)
    for every $k\in\{-2\gamma_j,1-2\gamma_j,2-2\gamma_j,2\gamma_j\}\cup\{-1,0\}$ and $m\in \{0,1\}$, 
    under hypothesis $\calH_{\Lk}$.
\end{enumerate}
\end{lemma}

The proof of Lemma~\ref{Convergence_Tps_long} is provided in Section~\ref{proof:convergence_tps_long}. It is quite technical so we give the main ideas and tools here. 
It relies on two auxiliary results for T-CKLS processes: Propositions~\ref{Preuve_holder}-\ref{Hitting_time}, whose proofs are postponed to Sections~\ref{proof:Preuve_holder}-\ref{proof:Hitting_time} respectively.
Proposition~\ref{Preuve_holder} is a property very commonly used in statistics for diffusion processes as well as the consequence of It\^o-formula in Proposition~\ref{prop:Ito}. 

\begin{prop}\label{Preuve_holder}
Assume that the T-CKLS $X$ is stationary.
Let $j\in\{0,\ldots,d\}$ and $m\geq 1$. Assume that the $m$-th moment of $\mu$ is finite. Then there exists a constant $C\in (0,\infty)$ such that for all $0\leq s<t<s+1$ it holds
$
    \Emu\left[|X_t-X_s|^m\right]\leq C(t-s)^{\sfrac{m}{2}}.
$
\end{prop}

\begin{prop}\label{prop:Ito}
Assume that the T-CKLS $X$ is stationary.
Let $j\in\{0,\ldots,d\}$ and $k\neq 0,1$.  
Assume that $\gamma_0\in [1/2,1]$ and that the $k,(k-1)$-th moments of $\mu$ are finite, as well as the moments $k-2(1-\gamma_i)$ and $2(k-1+\gamma_i)$-th moment for $\mu$ on $I_i$ with $i=0,d$. 
Then there exists a constant $C\in (0,\infty)$ such that for all $0\leq s<t<s+1$ it holds
$
    \Emu\left[|X_t^k-X_s^k| \mathds{1}_{I_j}(X_t) \right]\leq C(t-s)^{\sfrac{1}{2}}.
$
\\
Moreover, if $\mu$ admits finite $k,k-1,k-2(1-\gamma_j), 2(k-1+\gamma_j)$-th moment on $I_j$, then we have
$
    \Emu\left[|X_t^k-X_s^k| \mathds{1}_{I_j}((X_u)_{u\in [s,t]}) \right]\leq C(t-s)^{\sfrac{1}{2}}.
$
\end{prop}

Proposition~\ref{Hitting_time}, instead, is the remedy to the lack of knowledge of the finite dimensional distributions of TDs. 
Indeed, the quantities $\calM^{j,\cdot}, \calQ^{j,\cdot}$ consider only observations taking values on $I_j$ together with their following observation. 
Hence, in the proof of Lemma~\ref{Convergence_Tps_long}, one needs to bound the probability that the process crossed a threshold between two consecutive observations.

\begin{prop}\label{Hitting_time}
 Assume that the T-CKLS $X$ is stationary.
 Let $j\in\{0,\ldots,d\}$, $k\geq 1$, and $m\in\mathbb{R}$, such that the $m$-th moment of $\mu$ is finite on the set $I_j$.
 Then there exists a constant $C\in (0,\infty)$ such that for all $0\leq s<t<s+1$, we have:
\begin{equation}
\Emu\left[|X_s|^m\mathbb{P}_{X_s}\left( \tau_{I_j}^{\xi^{(j)}}<t-s\right)\mathds{1}_{I_j}(X_s) \right]\leq C(t-s)^{\sfrac{k}{2}},
\end{equation}
where $\xi^{(j)}$ is a CKLS process with parameters $(a_j,b_j,\sigma_j,\gamma_j)$ starting at $X_s$ and driven by a Brownian motion independent of $\mathcal{F}_s$ (denoted by $B$ as well),
and $\tau_{I_j}^{\xi^{(j)}}$ is the first hitting time of the interval $I_j$ for the process $\xi^{(j)}$.
\end{prop}

\section{Numerical Experiments} \label{sec:numerics}

In this section, we implement the MLE and QMLE based on discrete observations on simulated and US interest rates data.

\subsection{Simulated Data}
In this section, we investigate the efficiency of our estimators on simulated data. 
We simulate the T-CKLS process combining known Euler-Type schemes on different intervals such as the 
scheme in~\cite[equation~(3)]{alfonsi2005discretization} when the process is a CIR or a drifted version of the scheme in~\cite{taguchi2019euler} when the diffusion coefficient is non-linear. 
More precisely we use the following scheme. Given $X_0 \in (0,\infty)$ and $(G_k)_{k\in \mathbb{N}}$ a sequence of i.i.d.~standard Gaussian random variables, next we set 
$\Xn_0=X_0$ and, we define for all $k\in \mathbb{N}$ 
\[
    \Xn_{{(k+1)}/n} := \left| \Xn_{{k}/n} + \frac{1}{n} \left( a(\Xn_{{k}/n}) - b(\Xn_{{k}/n}) \Xn_{{k}/n}\right)+ \frac1{\sqrt{n}}\sigma(\Xn_{{k}/n}) (\Xn_{{k}/n})^{\gamma(\Xn_{{k}/n})}  G_k \right|. 
\]
Further discussion about the most suitable numerical scheme for T-CKLS is beyond the purpose of this paper.
To estimate the parameters from the simulated data, we use the estimators from discrete observations in Section~\ref{sec:estimators:discrete}.
The implementation has been done using Matlab and the parameters are as in Table~\ref{Parameters_simu}.

\begin{table}[H]
\begin{subtable}{0.45\textwidth}
\begin{center}
\begin{tabular}{llll}
\hline
$a_0$ & $b_0$  & $\sigma_0$  & $\gamma_0$  \\ \hline
0.3 & 0.2 & 0.2 & 0.5 \\ 
\end{tabular}
\end{center}
\end{subtable}
\begin{subtable}{0.45\textwidth}
\begin{center}
\begin{tabular}{llll}
\hline
$a_1$ & $b_1$  & $\sigma_1$  & $\gamma_1$  \\ \hline
0 & 0 & 0.4 & 0 \\  
\end{tabular}
\end{center}
\end{subtable}
\\
\par\bigskip
\begin{subtable}{0.45\textwidth}
\begin{center}
\begin{tabular}{llll}
\hline
$a_2$ & $b_2$  & $\sigma_2$  & $\gamma_2$  \\ \hline
0.3 & 0.2 & 0.2 & 0.5\\
\end{tabular}
\end{center}
\end{subtable}
\begin{subtable}{0.45\textwidth}
\begin{center}
\begin{tabular}{llll}
\hline
$r_1$ & $r_2$   \\ \hline
1 & 1.5 \\
\end{tabular}
\end{center}
\end{subtable}
\caption{Simulations parameters.}
\label{Parameters_simu}
\end{table}

Firstly, we illustrate Theorem \ref{Theorem_discret_cons_tcl} for the drift parameters $\theta_\star=(a_i,b_i)_{i=0}^{d=2}$, diffusion parameters $\sigma_\star$, $\gamma$ and two thresholds $r_1,r_2$ given in Table \ref{Parameters_simu}. 
We consider a process, which follows a CIR dynamic close to 0 and far away from 0 and is a BM on an intermediate bounded interval.
We simulate $10^4$ trajectories of the T-CKLS with two threshold. The set of numerical parameters is $(T,N)=(10^3,10^6)$ with starting condition determined as follows. As the process is supposed to be stationary, we first simulate one trajectory starting from $X_0>0$ chosen arbitrarily, say $X_0=1$, with the scheme $\Xn$ given above with $n =  k N/T$ for some $k\in \mathbb{N}\setminus \{0\}$ (we took $k=1$). Then, we consider the final value of the latter trajectory as initial condition of the $10^4$ trajectories.

\begin{figure}[H]
\hspace*{-1cm}    
	\begin{subfigure}[b]{0.5\linewidth}
		\centering \includegraphics[scale=0.6]{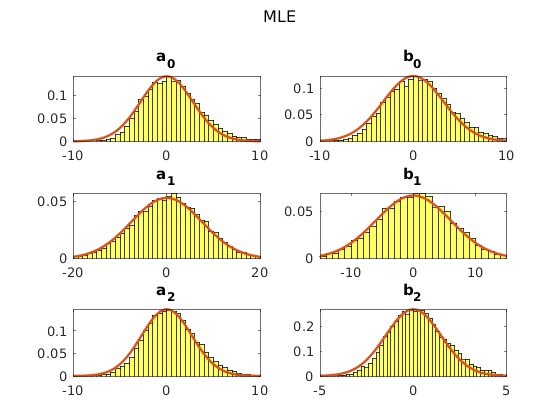}
	\end{subfigure}\hfill
	\begin{subfigure}[b]{0.5\linewidth}	
		\centering \includegraphics[scale=0.6]{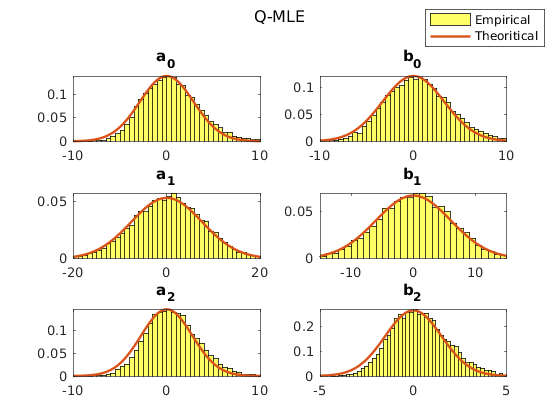}
	\end{subfigure}
	\caption{Asymptotic normality property in Theorem \ref{Theorem_discret_cons_tcl}, with parameters as in Table \ref{Parameters_simu}. We plot the theoritical distribution using \eqref{mesure_stationnaire} and compare with the empirical distribution on $10^4$ trajectories.}
	\label{CLT}
\end{figure}
Remark that, despite the fact that the set of numerical parameters $(T,N)$ does not satisfy the conditions for the asymptotic normality in Theorem~\ref{Theorem_discret_cons_tcl} ($T^2 =N$ instead of $T^2 \ll N$), numerics show good results. 

We compare the estimators by means of relative root-mean-square error (RMSE) and empirical bias (EB) in Table~\ref{tab:RMSE}.

\begin{table}[H]
\begin{center}
\begin{tabular}{clclcl}
\hline
Estimator & Relative RMSE & EB \\ \hline
  $(a_{T_N,N}^{0,(\Lk)},a_{T_N,N}^{1,(\Lk)},a _{T_N,N}^{2,(\Lk)} )$        & (0.3206, 0.2341, 0.2921)    & (0.0154, 0.0092, 0.0156) \\ \hline
$(b_{T_N,N}^{0,(\Lk)},b_{T_N,N}^{1,(\Lk)},b_{T_N,N}^{2,(\Lk)} )$         & (0.5500, 0.1857, 0.2439)    & (0.0173, 0.0069, 0.0093) \\ \hline
$(a_{T_N,N}^{0,(\qLk)},a_{T_N,N}^{1,(\qLk)},a_{T_N,N}^{2,(\qLk)} )$      & (0.3233, 0.2341, 0.2965)    & (0.0170, 0.0092, 0.0191) \\ \hline
  $(b_{T_N,N}^{0,(\qLk)},b_{T_N,N}^{1,(\qLk)},b_{T_N,N}^{2,(\qLk)} )$    & (0.5555,0.1857, 0.2478)    & (0.0182,0.0069, 0.0112) \\ \hline
$(\sigma_{T_N,N}^{0},\sigma_{T_N,N}^{1},\sigma_{T_N,N}^{2} )$        & (0.0088, 0.0087, 0.0015)    & (0.0087, -0.0034 , 0.001) \\ \hline
\end{tabular}
\caption{Table of the relative RMSE and EB for the estimator of $(\theta_\star,\sigma_\star)$ using the MLE \eqref{Estim_discret_multi_MLE}, QMLE \eqref{Estim_discret_multi_QMLE_1} and volatility estimator \eqref{estim_sigma_discret}.} \label{tab:RMSE}
\end{center}
\end{table}

Observe that the MLE gives a better estimation of the drift parameters, which is easily explained by the fact that the likelihood contains more information about the model. In general the MLE tends to have a better RMSE and EB than the QMLE. Applying the estimator on several data sets, we remark that the QMLE has a greater variance than the MLE.

\subsection{Interest rates analysis}\label{Interest rates analysis}

In this section, we apply our estimators to the ten year US treasury rate based on the Federal Reserve Bank's H15 data set. 
We exploit here the discretized version, considered in~\cite{Yu2020}, of the test for presence of thresholds of~\cite{su2017testing}. 
In the latter paper, one can find a computationally efficient approach for calibrating the p-value.
The test procedure, we describe below, simultaneously estimates and tests the threshold. 
Note that we use the drift MLE and $\sigma_\star$ estimators for CKLS process considered in this paper,
but we do not have results on threshold estimation, therefore we rely on the above mentioned results.
Indeed, while our Central Limit Theorem (CLT) in Theorem~\ref{Theorem_discret_cons_tcl} would allow us to construct a test to detect a threshold in the drift term (similarly to~\cite{mazzonetto2020drift}), additional theoretical results are required to estimate the new threshold.

\paragraph{General step.}
Let us describe the test to be applied at a general step and then the procedure describing the steps.
Suppose that there are $m$ thresholds, and that we look for the presence of an additional threshold on the $k$-th interval $I_k=(r_{k-1},r_{k})$ (we know the value of $r_{k-1},r_{k}$: either known or estimated in previous steps). We consider the hypothesis:
\begin{equation}
    \left\{\begin{array}{ll}
       H_0: \text{ Null hypothesis } \quad\quad\quad\quad m \text{ thresholds;}  \\
      H_1: \text{  Alternative hypothesis } \quad (m+1) \text{ thresholds.} 
    \end{array}\right.
    \label{Test}
\end{equation}  
Under the null hypothesis $H_0$, the model has $m$ thresholds. Under hypothesis $H_1$, there is an additional threshold $\bar{r}$ in the $k$-th regime, meaning that the sequence of thresholds becomes $-\infty=r_0<r_1<\ldots<r_{k-1}<\bar{r}<r_{k}<\ldots<r_m$. The quasi-likelihood ratio test statistic is given by
\[ 
    T=\sup_{\bar{r}\in[q_{20},q_{80}]} T (\bar{r})
\]
with
\begin{equation}
   T(\bar{r}):= 2\left(\qLk_{T_N,N}\left(\theta_{T_N,N}^{(H_1)}(m+1,\bar{r})\right)-\qLk_{T_N,N}\left(\theta_{T_N,N}^{(H_0)}(m)\right)\right)
    \label{test_stat}
\end{equation}
where $q_{20}$ and $q_{80}$ are $20$ and $80$ percentiles of the data in $[r_{k-1},r_k]$, the value 
 $\theta_{T_N,N}^{(H_1)}(m+1,\bar{r})$ is the drift MLE under 
 hypothesis $H_1$ with the additional threshold given by $\bar{r}$, and $\theta_{T_N,N}^{(H_0)}(m)$ is the drift MLE of the model under hypothesis $H_0$. The MLE is given in Proposition~\ref{prop:discrete_estimator}, with volatility parameter $\sigma_\star$ estimated by~\eqref{estim_sigma_discret}.

We compute the statistics for some values of $\bar{r}$, say $\bar{r}_j:= q_{20} (1-j/n) +  q_{80} j/n$, $j\in \{0,1,\ldots,n\}$ (we choose $n=10^3$). Next, we take as an estimator for the threshold $r_j$ the $\bar r_j$ which maximises $T(\bar{r})$ and the observed test statistics $T_{data}$ is then the quantity $T(\bar r_j)$.

The distribution of the test statistic \eqref{test_stat} is obtained using a bootstrap method. 
To compute the p-value we simulate $10^3$ trajectories of the process with the parameters under $H_0$, we compute $T_j$ the statistics on the $j$-th trajectory. Then the p-value is given by ${\#\{j\colon T_{data}<T_j\}}/{10^3}$.

We fix the significance level at the conventional $5\%$.

\paragraph{Procedure.}
We apply the test above in a sequential procedure, as follows. 
We first test for the presence of a threshold on the data: applying this test for $m=0$. 
If the test is significant, then we take as an estimator for the threshold $\widehat r_1$ the $\bar r$ which realises the maximum in the test statistics. This threshold divides the state space into two intervals. 
We then test the presence of thresholds on each of the two intervals, starting from the left to the right. On each interval, if the test is significant, we keep dividing the interval into two sub-intervals and so on. Once we do not have evidence of new thresholds in the interval we are considering, we go to the next interval.%

\paragraph{Application to ten year US treasury rate.}
We consider the ten year US Treasury rate and we adopt the convention that the daily time intervals is $\vd t=0.046$, where one unit of time represents one month. We assume that the data follow a T-CIR dynamics, \emph{i.e.}~$\gamma\equiv1/2$. 

\begin{figure}[H]
    \hspace{-0.5cm}
	\begin{subfigure}[b]{0.5\linewidth}
		\centering \includegraphics[scale=0.4]{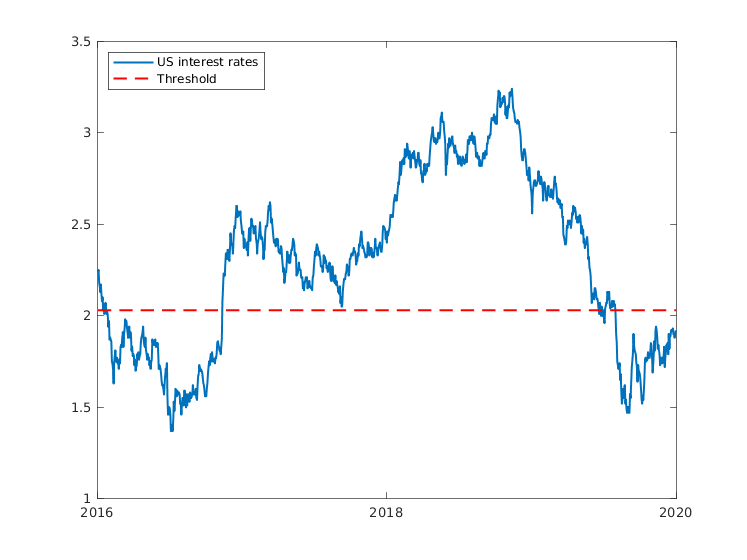}
	\end{subfigure}\hfill
	\begin{subfigure}[b]{0.5\linewidth}	
		\centering \includegraphics[scale=0.4]{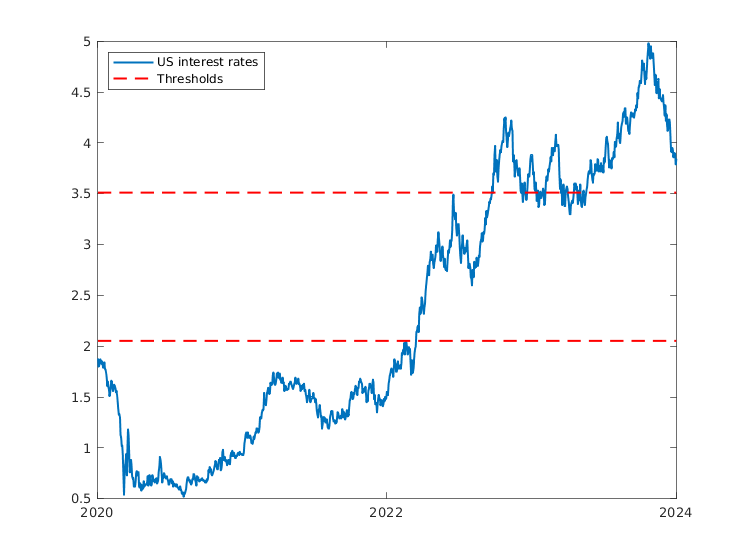}
	\end{subfigure}
	\caption{The figure shows the interest rate daily data (solid line) for the time window Jan 2016 - Dec 2019, and Jan 2020 - Jan 2024. The fitted thresholds are represented by the dashed lines.}
 \label{Test_fig}
\end{figure}

We consider the ten year US Treasury rate for two different time window: Jan 2016 - Dec 2019, and Jan 2020 - Jan 2024 represented in Figure~\ref{Test_fig}.

Let us consider the time window, Jan 2016 - Dec 2019.
The threshold test~\eqref{Test} for $m=0$ is significant and the threshold estimation is 
$r_1=2.0303$. 
We apply the test~\eqref{Test} to detect a threshold on $(0,r_1)$. It is not significant. The same conclusion holds testing for threshold presence on $(r_1,+\infty)$.

On the time window, Jan 2020 - Jan 2024, the threshold test~\eqref{Test} for $m=0$ is significant and we estimate the threshold $r_1=2.0507$.
There is no evidence of additional thresholds on $(0,r_1)$, and the null hypothesis is rejected for the existence of a threshold $r_2= 3.5112$ on $(r_1,+\infty)$. 
Instead, for the tests~\eqref{Test} with $m=2$ for finding a threshold on $(r_1,r_2)$ and $(r_2,+\infty)$, the null hypothesis $H_0$ is not-rejected.

\begin{table}[H]
\begin{subtable}{0.43\textwidth}
\subcaption*{Jan 2016 - Dec 2019}
\begin{tabular}{cc}
\hline
Estimator &  Value  \\ \hline
  $(a_{T_N,N}^{0,(\Lk)},a _{T_N,N}^{,(\Lk)} )$        & (1.6434 , 0.1713)     \\ \hline
$(b_{T_N,N}^{0,(\Lk)},b_{T_N,N}^{1,(\Lk)} )$         & (0.9410,  0.0723)     \\ \hline
$(\sigma_{T_N,N}^{0},\sigma_{T_N,N}^{1} )$        & (0.1616,  0.1053)     \\ \hline
$r_{T_N,N}^{1}$        & $\quad$$\quad$(2.0303)     \\ \hline
\end{tabular}
\end{subtable}
\begin{subtable}{0.4\textwidth}
\subcaption*{$\qquad\qquad\qquad$ Jan 2020 - Jan 2024} 
\begin{tabular}{cc}
\hline
Estimator & $\quad$$\quad$$\quad$$\quad$ Value  \\ \hline
  $(a_{T_N,N}^{0,(\Lk)},a_{T_N,N}^{1,(\Lk)},a _{T_N,N}^{2,(\Lk)} )$        & (0.2013,  0.5826, -0.0207)     \\ \hline
$(b_{T_N,N}^{0,(\Lk)},b_{T_N,N}^{1,(\Lk)},b_{T_N,N}^{2,(\Lk)} )$         & ( 0.1556, 0.0670, 0.0236)     \\ \hline
$(\sigma_{T_N,N}^{0},\sigma_{T_N,N}^{1},\sigma_{T_N,N}^{2} )$        & ( 0.2129,  0.2091, 0.1807)     \\ \hline
$(r_{T_N,N}^{1},r_{T_N,N}^{2})$        & $\qquad$(2.0507, 3.5112)     \\ \hline
\end{tabular}
\end{subtable}
\caption{Estimated parameters corresponding to Figure \ref{Test_fig}.}
\label{table_estim}
\end{table}

Therefore, we conclude that there is a single threshold in the time window, Jan 2016 - Dec 2019, and two thresholds in Jan 2020 - Jan 2024.
In Table \ref{table_estim}, we summarize the values obtained for each of the fitted parameters using the estimators.

\appendix

\section{Appendix: Properties of the process}\label{appendix}

Note that T-CLKS, solution to~\eqref{Multi_CKLS}, shows several behaviors; it may behave as an OU process on some intervals, a CIR or CKLS in others. 

The state space of T-CKLS is determined by the process behavior around 0, in particular at $I_0$ where it behaves as a standard CKLS process. The state space and the behavior at $0$ is provided in Lemma~\ref{lem:state_space}. As for standard CKLS process, if $\gamma_0\geq 1/2$, it is non-negative.

The regime of the process (transient, recurrent, positive-recurrent) is also determined by the behavior at $I_0$ and $I_d$ and so by the coefficients $a_0,b_0,\sigma_0,\gamma_0$ and $a_d,b_d,\sigma_d,\gamma_d$. We provide in Table~\ref{table_ergo} the conditions on the parameters such that the process is ergodic and we study the finiteness of the moments of the stationary distribution $\mu$ in Proposition~\ref{Moment_sup}.

Finally, we provide some auxiliary results on well posedness of some integrals of the T-CKLS appearing in the likelihood in Proposition~\ref{Proposition_integral}.

\subsection{State space and stationary distribution} \label{sec:property_stationary}

Ergodic T-CKLS process is a regular one-dimensional diffusion, so it is characterized by scale function $S$ and speed measure $m(x)\vd x$. The interested reader could refer to~\cite[II.4]{borodin2015handbook} for a summary or find more details, e.g.,~in~\cite[Chapter VII, Section 3]{revuz2013continuous}.  
The scale function is continuous, unique up to a multiplicative constant, and its derivative satisfies $S'(x)=\exp\left(-\int_{r_1}^x\frac{2(a(y)-b(y)y)}{\sigma(y)^2y^{2\gamma(y)}}\vd y\right)$. 
The speed measure is given by $m(x)\vd x = \frac{2}{(\sigma(x))^2 |x|^{2\gamma(x)} S'(x)}\vd x$.

The state of space of the T-CKLS process, denoted $I=\cup_{j=0}^d I_j$,  depends on the value of the parameters in $I_0$ ($a_0$, $\sigma_0$ and $\gamma_0$).  
 \begin{lemma} \label{lem:state_space}
     Let $X$ be the solution to the SDE~\eqref{Multi_CKLS}. 
     \begin{itemize}
        \item  If $\gamma_0=0$,  the state of space of the process is $I=\mathbb{R}$.
         \item 
         If $\gamma_0={1}/{2}$ and $0<a_0< {\sigma_0^2}/{2}$, the state space is $I=[0,+\infty)$ and the point ${0}$ is instantaneously reflecting.
         \item  If $\gamma_0\in ({1}/{2},1]$ or if $\gamma_0={1}/{2}$ and $a_0\geq {\sigma_0^2}/{2}$, then 
         $I=(0,+\infty)$ and ${0}$ is an unattainable boundary.
     \end{itemize}
 \end{lemma}

The regime of the process can be obtained by properties on the scale function and the spead measure (see, e.g.,~\cite[Exercice 3.15 in Chapter X]{revuz2013continuous} and \cite[Theorem 20.15]{kallenberg1997foundations}). The fact that $m(x) \vd x$ is a finite measure corresponds to positive recurrence (\emph{ergodicity}) of the process and 
the stationary measure is: 

\begin{equation}
    \label{mesure_stationnaire}
    \mu(\!\vd x)=\frac{m(x)}{\int_Im(y)\vd y}\vd x.
\end{equation}
The recurrent positivity property of the process only depends on the parameters below the first threshold (on $I_0$) and above the last threshold (on $I_d$). In the following table, we give conditions on the parameters $a_0,a_d$, $b_0,b_d$ and $\sigma_0,\sigma_d$ depending on the value of $\gamma_0,\gamma_d$ such that the process is ergodic (admits a stationary distribution).

\begin{table}[H]
\begin{subtable}{0.48\textwidth}
\begin{tabular}{|l|c|}
\hline
$\gamma_0=0$ &  $a_0\in \mathbb{R}$ and $b_0>0$, or \\
& $a_0>0$ and $b_0=0$ \\[1.9pt] \hline
$\gamma_0 \in [1/2,1)$ & $a_0>0$ and $ b_0 \in \mathbb{R}$
\\[1.9pt] \hline
$\gamma_0=1$ & $a_0>0$ and $ b_0 \in \mathbb{R}$, or
\\ & $a_0=0$ and $b_0< -{\sigma_0^2}/{2}$  \\[1.9pt] \hline
\end{tabular}
\end{subtable}
\begin{subtable}{0.52 \textwidth}
\begin{tabular}{|l|c|}
\hline
& $a_d \in \mathbb{R}$ and $b_d>0$, or
\\
\hline
$\gamma_d\in [0,1/2]$ &  $ a_d <0$ and $b_d=0$ \\[1.9pt] \hline
$\gamma_d\in(1/2,1)$ & $a_d\in \mathbb{R}$ and $b_d=0$ \\[1.9pt] \hline
$\gamma_d=1$ & $a_d \in \mathbb{R}$ and $b_d \in \left(-{\sigma_d^2}/{2} , 0\right]$ \\[1.9pt] \hline
\end{tabular}
\end{subtable}
\caption{Parameter conditions for ergodicity of T-CKLS process $X$ solution to~\eqref{Multi_CKLS}. 
}
\label{table_ergo}
\end{table}

The following proposition describes the behavior of the moments from the stationary distribution of the T-CKLS process in the ergodic regime. 

\begin{prop}\label{Moment_sup}
    Let $m\in \Rp$ and assume that the conditions in Table~\ref{table_ergo} (ensuring that $\mu$ is the stationary distribution) hold.
Then
$\mu$ admits finite $m$-th moment unless it holds simultaneously $b_d= 0$ and $\gamma_d\in [1/2,1]$ in which case the $m$-th moment is finite if
\begin{itemize}
    \item $\gamma_d=1/2$, $b_d=0$ and $a_d< -m\, {\sigma_d^2}/{2}$.
    \item $\gamma_d\in (1/2,1)$, $b_d=0$ and $m<\gamma_d-1/2<1$.
    \item $\gamma_d=1$, $b_d=0$ and $m\leq 1$. 
\end{itemize}
The measure $\mu$ admits finite $-m$-th moment unless $\gamma_0\in \{0,1/2\}$ in which case the $(-m)$-th moment is finite if 
\begin{itemize}
    \item $\gamma_0=1/2$ and $a_0> m \,  {\sigma_0^2}/{2}$.
    \item $\gamma_0=0$ and $m<1$. 
\end{itemize}
\end{prop}

\subsection{Well definiteness of the likelihood function}

The following proposition describes the behavior of various integrals of the T-CKLS process. It establishes whether the likelihood~\eqref{Likelihood} is well defined or not. Properties of this kind have also been considered in~\cite{mishura2024low}.

\begin{prop} \label{Proposition_integral}
    Let $X$ be solution to the SDE \eqref{Multi_CKLS}. 
    \begin{enumerate}[(a)]
        \item If $\gamma_0={1}/{2}$ and $a_0\geq {\sigma_0^2}/{2}$ or if $\gamma_0\in\left({1}/{2},1\right]\cup\{0\}$, then 
        \begin{equation}
           \forall t\geq 0, \quad \int_0^t\frac{1}{X_s^{2\gamma_0}}\mathds{1}_{I_0}(X_s)\vd s <\infty \quad \mathbb{P}_{x_0}-\text{ a.s.}.
        \end{equation}
        \item If $\gamma_0={1}/{2}$ and $a_0<{\sigma_0^2}/{2}$, then 
        \begin{equation}
            \forall t\geq 0, \quad \mathbb{P}_{x_0}\left(\int_0^t\frac{1}{X_s^{2\gamma_0}}\mathds{1}_{I_0}(X_s)\vd s=\infty\right)>0.
        \end{equation}
    \end{enumerate}
\end{prop}
\begin{proof}
    If $\gamma_0=0$, it is trivial. 
    When $\gamma_0\neq 0$, the first item follows from the fact that 0 is an unattainable boundary and continuity of the trajectories: the image of $[0,t]$ through each trajectory $s \mapsto X_s(\omega)$ is a compact of $ ]0,\infty[$.
    The second item, for $\gamma_0=1/2$ has been proven in~\cite{ben2013asymptotic} by using properties of the Laplace transform.
\end{proof}

\section{Appendix: Proofs}\label{Proof}
In this section, we first derive $\frM_{T_N,N}^{.,m}$ in~\eqref{frak_M_1_discr} and then we prove Theorem~\ref{Thm_sigma}, Theorem~\ref{Theorem_discret_cons_tcl}, Lemma~\ref{Convergence_Tps_long}, Proposition~\ref{Preuve_holder} and Proposition~\ref{Hitting_time}.

\subsection{Likelihood discretization: $\calM_{T_N,N}^{.,m}$ versus $\frM_{T_N,N}^{.,m}$} 

We derive $\frM_{T_N,N}^{.,m}$ in~\eqref{frak_M_1_discr} as a discretisation of an expression $\mathbb{P}$-a.s.~equal to $\calM^{j,m}_T$, which would not involve any term $\calM^{j,k}$ except if $k=0$. The notation, in particular the functions $f_{j,m}$, $\frak{f}_j$, have been introduced together with equation~\eqref{frak_M_1_discr}.

\begin{lemma}\label{lemma_disc}
Let $T\in (0,\infty)$, $j\in\{1\ldots,d\}$ and $m\in\{-2\gamma_j,1-2\gamma_j\}\setminus \{0\}$. It holds $\mathbb{P}$-a.s.~ that
\begin{align}
    \calM^{0,m}_{T}= & f_{0,m+1}(X_0)-f_{0,m+1}(X_T)-\frac{m}{2}\sigma_0^2\calQ_T^{0,m+2\gamma_0-1}   + r_1^m\left(\calM_{T}^{0,0} +\frak{f}_0 \right),
    \label{frak_M_2_0}
\end{align}
\begin{align}
    \calM^{d,m}_{T}= & f_{d,m+1}(X_T)-f_{d,m+1}(X_0)-\frac{m}{2}\sigma_d^2\calQ_T^{d,m+2\gamma_d-1}   + r_{d}^m\left(\calM_{T}^{d,0} - \frak{f}_d \right),
    \label{frak_M_2_d}
\end{align}
and for $j\in\{1,\ldots, d-1\}$:
\begin{align} 
    \calM^{j,m}_{T} &= f_{j,m+1}(X_T) - f_{j,m+1} (X_0)-\frac{m}{2}\sigma_j^2\calQ_{T}^{j,m+2\gamma_j-1}\nonumber \\
    & \quad  + r_j^m \calM_{T}^{j,0} + r_{j+1}^m \frak{f}_{j+1}- r_{j}^m \frak{f}_{j} - (r_{j+1}^m - r_j^m) \sum_{k=j+1}^{d} \calM_{T}^{k,0}.
    \label{frak_M_1}
\end{align}

\end{lemma}
\begin{proof}
 We prove only the case $j\in\{1,\dots,d-1\}$, and the case $j\in \{0,d\}$ works similarly.
 Applying It\^o-Tanaka formula yields the a.s.~equality 
 \begin{equation}
 f_{j,m+1}(X_T)=f_{j,m+1} (X_0)+\frac{m}{2}\sigma_j^2\calQ^{j,m+2\gamma_j-1}+ \calM^{j,m}_{T}+\frac12\left( r_j^m L_T^{r_j}(X)- r_{j+1}^m L_T^{r_{j+1}}(X)\right).
 \label{eq:discret_M}
 \end{equation}

 Applying a second time It\^o-Tanaka formula to $\max{(X_T,r_j)}$ shows that a.s.~$
       L^{r_j}_T(X)  
        = 2 \frak{f}_j - 2 \sum_{i=j}^d \calM_T^{i,0}.
     \label{eq:discret_M_loc1}
 $
Combining this with \eqref{eq:discret_M} completes the proof.
\end{proof}

To obtain $\frM^{j,m}_{T_N,N}$, we just consider the discretized versions of the right-hand side of~\eqref{frak_M_2_0}, \eqref{frak_M_2_d} and \eqref{frak_M_1} in Lemma~\ref{lemma_disc}: we replace the quantities $\calM^{j,0}$ and $\calQ^{j,m+2\gamma_j -1}$ by their discretized versions~\eqref{Quantity_likelihood_discrete_multi}.

\subsection{Proof of Theorem \ref{Thm_sigma}}

We study the asymptotic behavior of $(\sigma_{T_N,N})^2-\sigma^2_\star$. For $j\in\{0,\ldots,d\}$, by Proposition \ref{Estimateur_sigma_cont}, on the event $\{\calQ_T^{j,0}>0\}$, we have:
\begin{equation}
(\sigma^j_{T_N,N})^2-\sigma^2_j =\sigma_j^2 \frac{\calQ_{T_N}^{j,2\gamma_j}-\calQ_{T_N,N}^{j,2\gamma_j}}{\calQ_{T_N,N}^{j,2\gamma_j}}+\frac{\frQ^j_{T_N,N}-\frQ^j_{T_N}}{\calQ_{T_N,N}^{j,2\gamma_j}}.
\end{equation}
 By equations~\eqref{discret_bracket2}-\eqref{discret_bracket1}, we have:
 $$\left|\frQ^j_{T_N}-\frQ^j_{T_N,N}\right|\leq 2\left|\calM_{T_N}^{j,1}-\calM_{T_N,N}^{j,1}\right| +2\max(|r_j|,|r_{j+1}|)\sum_{i=j}^d \left|\calM_{T_N}^{i,0}-\calM_{T_N,N}^{i,0}\right|,$$
 for $j\in\{1,\ldots,d-1\}$ and for $j\in\{0,d\}$
  $$\left|\frQ^j_{T_N}-\frQ^j_{T_N,N}\right|\leq 2\left|\calM_{T_N}^{j,1}-\calM_{T_N,N}^{j,1}\right| +2|r_{(j+1)\wedge d}| \left|\calM_{T_N}^{j,0}-\calM_{T_N,N}^{j,0}\right|.$$
  Then, we conclude by making use of item \ref{item_0} in Lemma \ref{Convergence_Tps_long} and the fact that\\ $\mathbb{P}\left(\lim_{T\to+\infty}\calQ_T^{j,0}>0\right)=1$.

\subsection{Proof of Theorem \ref{Theorem_discret_cons_tcl}} 
For all $N\in\mathbb{N}$ and $\ell\in\{\Lk,\qLk\}$ 
it holds that
\[
    \left(\theta^{(\ell)}_{T_N,N}-\theta_{\star}\right)=\left(\theta^{(\ell)}_{T_N,N} -\theta_{T_N}^{(\ell)}  \right)+ \left(\theta^{(\ell)}_{T_N} -\theta_{\star} \right) .
\]
The second term, on the right-hand side of the equality, provides the asymptotic behavior by applying Theorem \ref{Theorem_cont}. We now show that the first vanishes. Before that let us just mention that the denominator of all these quantities a.s.~never vanish (by Cauchy-Schwarz inequality). 
In the case $\ell=\Lk$, for $j\in\{0,\ldots,d\}$, using equations \eqref{estimateurs_likelihood_continue_multi} and \eqref{Estim_discret_multi_MLE},  each component of the first term of the right-hand side of the last equation can be rewritten as follows:
\[\begin{split}
    & \left(\frac{\calQ_{T_N,N}^{j,k} }{\calQ_{T_N,N}^{j,-2\gamma_j}\calQ_{T_N,N}^{j,2-2\gamma_j}-(\calQ_{T_N,N}^{j,1-2\gamma_j})^2  }-\frac{\calQ_{T_N}^{j,k} }{\calQ_{T_N}^{j,-2\gamma_j}\calQ_{T_N}^{j,2-2\gamma_j}-(\calQ_{T_N}^{j,1-2\gamma_j})^2  }\right)\frM_{T_N}^{j,m}
    \\
    & \qquad \qquad \qquad +\frac{\calQ_{T_N,N}^{j,k}( \frM_{T_N,N}^{j,m}-\frM_{T_N}^{j,m})}{\calQ_{T_N,N}^{j,-2\gamma_j}\calQ_{T_N,N}^{j,2-2\gamma_j}-(\calQ_{T_N,N}^{j,1-2\gamma_j})^2},
\end{split}\]
with $k\in\{-2\gamma_j,1-2\gamma_j,2-2\gamma_j\}$ and $m\in\{-2\gamma_j,1-2\gamma_j\}$. 
Let $\bar{r}:=\max\{|r_j| \colon j=1,\ldots, d\}$. 
Then, using formula \eqref{frak_M_2_0}, \eqref{frak_M_2_d} and \eqref{frak_M_1}, for $j\in\{1,\ldots,d-1\}$, we have
\begin{align}
    |\frM_{T_N}^{j,m}-\frM_{T_N,N}^{j,m}|&\leq \bar{r}^{m}\sum_{i=j}^d|\calM_{T_N}^{i,0}-\calM_{T_N,N}^{i,0}|
    + \frac{|m|}{2}\left| \left(\sigma_{T_N,N}^j\right)^2-\sigma_j^2\right| \left|\calQ_{T_N}^{j,m+2\gamma_j-1}-\calQ_{T_N,N}^{j,m+2\gamma_j-1}\right| ,
\end{align}
and for $j\in\{0,d\}$:
\begin{align}
   |\frM_{T_N}^{j,m}-\frM_{T_N,N}^{j,m}|
   &\leq \bar{r}^{m} |\calM_{T_N}^{j,0}-\calM_{T_N,N}^{j,0}|+\frac{ |m|}{2} \left| \left(\sigma_{T_N,N}^j\right)^2-\sigma_j^2\right|\left|\calQ_{T_N}^{j,m+2\gamma_j-1}-\calQ_{T_N,N}^{j,m+2\gamma_j-1}\right| .
\end{align}
This, Item~\ref{item_2} in Lemma~\ref{Convergence_Tps_long}, and Theorem \ref{Thm_sigma} ensure that 
\[
    T_N^{-{1}/{\lambda}} \left|\calQ_{T_N,N}^{j,k}-\calQ_{T_N}^{j,k}\right|\xrightarrow[N\to+\infty]{\mathbb{P}} 0
    \quad \text{and} \quad 
    T_N^{-{1}/{\lambda}} \left|\frM_{T_N,N}^{j,m}-\frM_{T_N}^{j,m}\right|\xrightarrow[N\to+\infty]{\mathbb{P}} 0
\]
with $\lambda=1$ to get the consistency of the MLE and $\lambda=2$ for the speed of convergence. 
In the case $\ell=\qLk$, the first term works analogously, exploiting Item~\ref{item_1} in Lemma~\ref{Convergence_Tps_long}. 

To establish the LAN property, we proceed with a similar reasoning. Let $\Delta \theta_\star:=(\Delta a_j,\Delta b_j)_{j=0}^d$. Then, we decompose the discretized log-likelihood ratio as 
\begin{align}
  & \ln\frac{ \Lk_{T_N,N}(\theta_\star+\frac{1}{\sqrt{T_N}}\Delta\theta_\star; \sigma_{T_N,N}, \gamma)}{\Lk_{T_N,N}(\theta_\star; \sigma_{T_N,N}, \gamma)}
  \\ & \quad = \ln\frac{ \Lk_{T_N,N}(\theta_\star+\frac{1}{\sqrt{T_N}}\Delta\theta_\star; \sigma_{T_N,N}, \gamma)}{\Lk_{T_N,N}(\theta_\star; \sigma_{T_N,N}, \gamma)}-\ln\frac{ \Lk_{T_N}(\theta_\star+\frac{1}{\sqrt{T_N}}\Delta\theta_\star; \sigma, \gamma)}{\Lk_{T_N}(\theta_\star; \sigma, \gamma)}
  \\ & \qquad +\ln\frac{ \Lk_{T_N}(\theta_\star+\frac{1}{\sqrt{T_N}}\Delta\theta_\star; \sigma, \gamma)}{\Lk_{T_N}(\theta_\star; \sigma, \gamma)}.
\end{align}
Similarly, for discretized the quasi-likelihood, we write
\begin{align}
    & \qLk_{T_N,N}\left(\theta_\star+\frac{1}{\sqrt{T_N}}\Delta\theta_\star\right)-\qLk_{T_N,N}(\theta_\star)
    \\ 
    & \quad = \qLk_{T_N,N}\left(\theta_\star+\frac{1}{\sqrt{T_N}}\Delta\theta_\star\right)-\qLk_{T_N}\left(\theta_\star+\frac{1}{\sqrt{T_N}}\Delta\theta_\star\right)+\qLk_{T_N}(\theta_\star) -\qLk_{T_N,N}(\theta_\star)\\
    & \qquad +\qLk_{T_N}\left(\theta_\star+\frac{1}{\sqrt{T_N}}\Delta\theta_\star\right)-\qLk_{T_N}(\theta_\star).
\end{align}

Applying Theorem~\ref{Theorem_cont}, Theorem~\ref{Thm_sigma}, together with Lemma~\ref{Convergence_Tps_long} completes the proof.

\subsection{Proof of Lemma~\ref{Convergence_Tps_long}} \label{proof:convergence_tps_long}

Let us introduce the round ground notation $\lfloor t\rfloor_{\Delta_N}:=t_i$ for $t\in[t_i,t_{i+1})\subseteq [t_i,t_i+\Delta_N]$. Moreover, without loss of generality, we assume $T_N\leq N$ and $\Delta_N \leq 1$ for all $N\in\mathbb{N}$. 

Analogously to~\cite{mazzonetto2020drift}, 
the proof relies on Propositions-\ref{Preuve_holder}-\ref{Hitting_time}. 
Basically, for $\lambda\in\{1,2\}$, the proof of Lemma~\ref{Convergence_Tps_long} reduces to prove that the following integrals are $o(T_N^{1/\lambda})$
\begin{equation}\label{1_condition}
    \int_0^{T_N} \Emu\left[|X_t-X_{\lfloor t\rfloor_{\Delta_N}}|^p \right]\vd t, \qquad 
    \int_0^{T_N} \Emu\left[|X_t^k-X_{\lfloor t\rfloor_{\Delta_N}}^k| \mathds{1}_{I_j}(X_t) \right]\vd t,
    \quad \text{and}
\end{equation}
\begin{equation}\label{eq:general_sans_tps_att}
    \int_0^{T_N} \Emu\left[|X_{\lfloor t\rfloor_{\Delta_N}}|^q\mathds{1}_{X_t\notin I_j, X_{\lfloor t\rfloor_{\Delta_N}}\in I_j } \right]\vd t     \quad \text{and} \quad \int_0^{T_N} \Emu\left[|X_{\lfloor t\rfloor_{\Delta_N}}|^q \mathds{1}_{X_t\in I_j, X_{\lfloor t\rfloor_{\Delta_N}}\notin I_j}  \right]\vd t
\end{equation}
for some suitable $p>0$ and $q\in \mathbb{R}$ and specific assumptions, all depending on the value of $m$ and $k$ appearing in $\calM^{j,m}$ and $\calQ^{j,k}$ in the different items of the statement we are proving. Let us use the notation $\mathbf{J}_{k,j}(t) := \Emu\left[|X_t^k-X_{\lfloor t\rfloor_{\Delta_N}}^k| \mathds{1}_{X_t \in I_j} \right]$ and 
\[\mathbf{Y}_{q,j}(t) :=\Emu\left[|X_{\lfloor t\rfloor_{\Delta_N}}|^q\mathbb{P}_{X_{\lfloor t\rfloor_{\Delta_N}}}\left( \tau_{I_j}^{\xi^{(j)}}<\Delta_N\right) \mathds{1}_{I_j}(X_{\lfloor t\rfloor_{\Delta_N}})\right]\]
where $\xi^{(j)}$ is a CKLS process with parameter $(a_j,b_j,\sigma_j,\gamma_j)$ starting at $\xi^{(j)}_0=X_{\lfloor t\rfloor_{\Delta_N}}$ but driven by a BM, (denoted by $B$ as well), independent of $\mathcal{F}_{\lfloor t\rfloor_{\Delta_N}}$
and $\tau_{I_j}^{\xi_j}$ is the first hitting time of the interval $I_j$ for the process $\xi^{(j)}$.

Let us consider equation \eqref{eq:general_sans_tps_att}. 
The tower property of conditional expectation and Markov property imply that
\begin{align}\label{eq:max:g}
 \Emu\left[|X_{\lfloor t\rfloor_{\Delta_N}}|^q \mathds{1}_{X_t\notin I_j, X_{\lfloor t\rfloor_{\Delta_N}}\in I_j}\right]&\leq \mathbf{Y}_{q,j}(t)
 \quad \text{and} \quad 
 \Emu\left[|X_{\lfloor t\rfloor_{\Delta_N}}|^q \mathds{1}_{X_t\in I_j, X_{\lfloor t\rfloor_{\Delta_N}}\notin I_j}\right] \leq \sum_{\substack{i=0\\ i\neq j}}^d \mathbf{Y}_{q.i}(t).
\end{align} 

Proposition~\ref{Hitting_time}, when it can be applied (i.e.~if $\mu$ admits finite moments of order $q$) ensures that 
\begin{equation} \label{eq:gsta_bound}
     \int_0^{T_N}  \mathbf{Y}_{q,j}(t)
     \vd t \leq C T_N \Delta_N
\end{equation}
for some positive constant $C$ independent of $N$.

The term $\mathbf{J}_{k,j}$ is mainly dealt with Proposition~\ref{prop:Ito} and we discuss it in detail in the next sections, where we see that sometimes it is convenient to bound it with terms involving $\mathbf{Y}_{k,\cdot}$, $\mathbf{Y}_{0,\cdot}$ and taking into account that, if $\mu$ admits finite $\max(1,p)$-th moment ($p>0$), 
\begin{equation} \label{eq:1_condition_bound}
    \int_0^{T_N} \Emu\left[|X_t-X_{\lfloor t\rfloor_{\Delta_N}}|^p \right]\vd t
    \leq 
   \int_0^{T_N} \Emu\left[|X_t-X_{\lfloor t\rfloor_{\Delta_N}}|^{\max(p,1)} \right]^{\min(p,1)}\vd t 
   \leq C T_N \Delta_N^{p/2}
\end{equation}
for some positive constant $C$ independent of $N$.

The assumption $\lim_{N\to +\infty} T_N^{\lambda-1} \Delta_N=0$ ensures that~\eqref{eq:1_condition_bound} and~\eqref{eq:gsta_bound} are $o(T^{1/\lambda})$.

The proof of each item of~Lemma~\ref{Convergence_Tps_long} is then reduced to determine the corresponding values of $p$ and $q$ in~\eqref{eq:1_condition_bound} and~\eqref{eq:gsta_bound}. By doing so, we deduce the assumptions $\calH_\Lk$ and $\calH_{\qLk}$ on moments of $\mu$.
We consider separately the two quantities in~\eqref{Conv_i_m} and then we summarize the items of Lemma~\ref{Convergence_Tps_long} in Section~\ref{sec:lemma52_proof}.

\subsubsection{Dealing with $\Emu\left[|\calM_{T_N}^{j,m}-\calM_{T_N,N}^{j,m}|\right]$} \label{sec:proof:lemme5_2:M}

In this section, we use~\eqref{eq:1_condition_bound}-\eqref{eq:gsta_bound} to obtain the assumptions on the moments of $\mu$ under which
$\Emu\left[|\calM_{T_N}^{j,m}-\calM_{T_N,N}^{j,m}|\right]$ is $o(T_N^{1/\lambda})$
for $m\in \{0,1\}$, $j\in \{0,\ldots,d\}$.

As for $\calQ^{j,k}$, we can rewrite 
$\calM_{T_N,N}^{j,m}-\calM_{T_N}^{j,m}$ as follows
\begin{align}
\calM_{T_N,N}^{j,m}-\calM_{T_N}^{j,m}&=\int_0^{T_N}\left(X_t^m\mathds{1}_{I_j}(X_t)-X^m_{\lfloor t\rfloor_{\Delta_N}}\mathds{1}_{I_j}(X_{\lfloor t\rfloor_{\Delta_N}}) \right)(a(X_t)-b(X_t)X_t)\vd t\\
&\qquad +\int_0^{T_N}\left(X_t^m\mathds{1}_{I_j}(X_t)-X^m_{\lfloor t\rfloor_{\Delta_N}}\mathds{1}_{I_j}(X_{\lfloor t\rfloor_{\Delta_N}}) \right)\sigma(X_t)(X_t)^{\gamma(X_t)}\vd B_t.
\end{align}

Using Triangular inequality, Hölder's inequality, and Itô-isometry, we obtain
\begin{align}
&\Emu\left[|\calM_{T_N}^{j,m}-\calM_{T_N,N}^{j,m}|\right]\leq \int_0^{T_N} \Emu\left[|X_t^m-X_{\lfloor t\rfloor_{\Delta_N}}^m|(|a_j|+b_j|X_{\lfloor t\rfloor_{\Delta_N}}|+b_j|X_t-X_{\lfloor t\rfloor_{\Delta_N}}|)\right]\vd t\nonumber\\
&+\int_0^{T_N}\Emu\left[|X_{\lfloor t\rfloor_{\Delta_N}}|^m\mathds{1}_{X_t\in I_j, X_{\lfloor t\rfloor_{\Delta_N}}\notin I_j}\left(\max_{i=0,\ldots,d} |a_i|+\max_{i=0,\ldots,d} |b_i|\left(|X_{\lfloor t\rfloor_{\Delta_N}}|+|X_t-X_{\lfloor t\rfloor_{\Delta_N}}|\right)\right) \right]\vd t \nonumber\\
&+\int_0^{T_N}\Emu\left[|X_{\lfloor t\rfloor_{\Delta_N}}|^m\mathds{1}_{X_t\notin I_j, X_{\lfloor t\rfloor_{\Delta_N}}\in I_j }\left(\max_{i=0,\ldots,d} |a_i|+\max_{i=0,\ldots,d} |b_i|\left(|X_{\lfloor t\rfloor_{\Delta_N}}|+|X_t-X_{\lfloor t\rfloor_{\Delta_N}}|\right)\right) \right]\vd t
\nonumber\\
&+\sqrt{2}\max_{i=0,\ldots,d}(\sigma_i)\left(\int_0^{T_N}\Emu\left[(X_t^m-X^m_{\lfloor t\rfloor_{\Delta_N}})^{2}(X_t)^{2\gamma(X_t)}+X_{\lfloor t\rfloor_{\Delta_N}}^{2m}(X_t)^{2\gamma(X_t)}  \mathds{1}_{X_t\in I_j, X_{\lfloor t\rfloor_{\Delta_N}}\notin I_j} \right.\right.\nonumber \\
&\left.\left.  +X_{\lfloor t\rfloor_{\Delta_N}}^{2m}(X_t)^{2\gamma(X_t)}\mathds{1}_{X_t\notin I_j, X_{\lfloor t\rfloor_{\Delta_N}}\in I_j }        \right]\vd t\right)^{1/2}.
\label{M}
\end{align}

\paragraph{The case $m=0$.} $\mu$ has finite $\max(1, 2\gamma_d)$-th moment. 

\noindent Inequality~\eqref{M} involves terms of the kind~\eqref{1_condition} with $p=1$ and~\eqref{eq:general_sans_tps_att} with $q\in \{0,1\}$ but also
\begin{align}
    & \left(\int_0^{T_N}\Emu\left[ (X_t)^{2\gamma(X_t)}  (\mathds{1}_{X_t\in I_j, X_{\lfloor t\rfloor_{\Delta_N}}\notin I_j} + \mathds{1}_{X_t\notin I_j, X_{\lfloor t\rfloor_{\Delta_N}}\in I_j })\right]\vd t\right)^{1/2}.
\end{align}
Since $I_j$ is not bounded only if $j=d$ or $j=0$ with $\gamma_0=0$, the above quantity is bounded from above for all $j\neq d$ and for $j=d$ if $\gamma_d=0$ by terms~\eqref{eq:general_sans_tps_att} with $q=0$.
If $j=d$ and $\gamma_d\neq 0$, its square is bounded by a sum of terms as~\eqref{1_condition} with $p=2\gamma_d$ and~\eqref{eq:general_sans_tps_att} with $q=2\gamma_d$:
\begin{equation}
    \Emu\left[ (X_t)^{2\gamma_d} \mathds{1}_{X_t \in I_d, X_{\lfloor t\rfloor_{\Delta_N}}\notin I_d } \right] 
    \leq    C \Emu\left[ |X_t-X_{\lfloor t \rfloor_{\Delta_N}}|^{2\gamma_d} + X_{\lfloor t \rfloor_{\Delta_N}}^{2\gamma_d} \mathds{1}_{X_t \in I_d, X_{\lfloor t\rfloor_{\Delta_N}}\notin I_d } \right]
    \end{equation}
for some positive constant $C$.
The latter inequality is derived from sub-additivity (if $\gamma_d \leq 1/2$) or Jensen's inequality (if $\gamma_d>1/2$).

Inequalities~\eqref{eq:1_condition_bound}-\eqref{eq:gsta_bound} ensure that
\begin{align}
    &\Emu\left[|\calM_{T_N}^{j,0}-\calM_{T_N,N}^{j,0}|\right]
    \leq C_1  T_N \sqrt{\Delta_N} + C_2 \sqrt{ T_N \Delta_N^{\max(1,2\gamma_d)/2} } \leq   (C_1+C_2)  T_N \sqrt{\Delta_N}
\end{align}
for some positive constants $C_1,C_2$.

\paragraph{The case $m=1$.} $\mu$ admits finite $2(1+\gamma_d)$-th moment. 

\noindent Inequality~\eqref{M} involves terms of the kind~\eqref{1_condition} with $p \in \{1,2\}$ and~\eqref{eq:general_sans_tps_att} with $q\in \{1,2\}$ but also the square root of the following term 
\begin{align}
    & 
    \int_0^{T_N}\Emu\left[ (X_t -X_{\lfloor t\rfloor_{\Delta_N}})^2 X_t^{2\gamma(X_t)}+ X_{\lfloor t\rfloor_{\Delta_N}}^2 X_t^{2\gamma(X_t)}(\mathds{1}_{X_t\in I_j, X_{\lfloor t\rfloor_{\Delta_N}}\notin I_j} + \mathds{1}_{X_t\notin I_j, X_{\lfloor t\rfloor_{\Delta_N}}\in I_j })\right]\vd t.
\end{align}
Similarly to the case $m=0$, by the fact that $X_s^{2\gamma_j} \mathds{1}_{I_j}(X_s)$ is not bounded from above (up to multiplicative constant) by $\mathds{1}_{I_j}(X_s)$ only if $j=d$, the fact that the same holds for $X_s^2 \mathds{1}_{I_j}(X_s)$ unless $j=d$ or $j=0$ and $\gamma_0=0$, and by sub-additivity of the square-root, we can reduce to 
\begin{align}
    &
    \left( 
    \int_0^{T_N} \Emu\left[ (X_t -X_{\lfloor t\rfloor_{\Delta_N}})^2 
    + \mathds{1}_{X_t\in I_j, X_{\lfloor t\rfloor_{\Delta_N}}\notin I_j} + \mathds{1}_{X_t\notin I_j, X_{\lfloor t\rfloor_{\Delta_N}}\in I_j }\right]
    \vd t \right)^{1/2}
    \\
    & \quad +
    \left(
    \int_0^{T_N} \Emu\left[
    (X_t -X_{\lfloor t\rfloor_{\Delta_N}})^{2} X_t^{2\gamma_d} 
    \mathds{1}_{I_d}(X_t)
    +
    X_{\lfloor t\rfloor_{\Delta_N}}^2 X_t^{2\gamma_d} \mathds{1}_{X_t\in I_d, X_{\lfloor t\rfloor_{\Delta_N}} \in  I_0}  \right] \right)^{1/2}
    \\
    & \quad + \left(
    \int_0^{T_N} \Emu\left[ 
   X_{\lfloor t\rfloor_{\Delta_N}}^2 \mathds{1}_{X_t\notin I_d, X_{\lfloor t\rfloor_{\Delta_N}}\in I_d }
   + 
   X_t^{2\gamma_d} \mathds{1}_{X_t\in I_d, X_{\lfloor t\rfloor_{\Delta_N}}\notin I_d \cup I_0}  
    \right] \right)^{1/2}.
\end{align}
The first line involves terms like~\eqref{1_condition} with $p=2$ and~\eqref{eq:general_sans_tps_att} with $q=0$. The last line shows a term~\eqref{eq:general_sans_tps_att} with $q=2$ and a term appearing for the case $m=0$.
The second line, by subadditivity and  H\"older's inequality ($(1+\gamma_d)^{-1}+ (1+\gamma_d)^{-1}\gamma_d=1$), 
is bounded from above by
\begin{align}
    & 
    \left(\int_0^{T_N} 
     \Emu\left[ (X_t -X_{\lfloor t\rfloor_{\Delta_N}})^{2(1+\gamma_d)}\right]^{1/(1+\gamma_d)} 
    \Emu\left[X_t^{2 (1+\gamma_d)}  
    \right]^{\gamma_d/(1+\gamma_d)} \vd t 
    \right)^{1/2}
    \\
    & \quad + \left(\int_0^{T_N} 
      \Emu\left[ X_{\lfloor t\rfloor_{\Delta_N}}^{2(1+\gamma_d)} \mathds{1}_{X_t \in I_d, X_{\lfloor t\rfloor_{\Delta_N}}\in I_0 } \right]^{1/(1+\gamma_d)} 
        \Emu\left[X_t^{2 (1+\gamma_d)} \right]^{\gamma_d/(1+\gamma_d)} \vd t 
    \right)^{1/2}
\end{align}
where we recognise~\eqref{1_condition} and~\eqref{eq:general_sans_tps_att} with $p=q=2(1+\gamma_d)$.

By~\eqref{eq:1_condition_bound} and~\eqref{eq:gsta_bound} we deduce
$\Emu\left[|\calM_{T_N}^{j,1}-\calM_{T_N,N}^{j,1}|\right]$ is bounded from above by 
\begin{equation}
    C_1 ( T_N \sqrt{\Delta_N} + T_N \Delta_N + \sqrt{T_N \Delta_N} + \sqrt{T_N \Delta_N^{1/2(1+\gamma_j)}} ) \leq 4 C_1 T_N \sqrt{\Delta_N} 
\end{equation}
for some positive constant $C_1$.

\subsubsection{Dealing with 
$\Emu\left[|\calQ_{T_N}^{j,k}-\calQ_{T_N,N}^{j,k}|\right]$.}

In this section, using~\eqref{eq:1_condition_bound}-\eqref{eq:gsta_bound}, we precise under which assumptions on the moments of $\mu$ it holds that
$\Emu\left[|\calQ_{T_N}^{j,k}-\calQ_{T_N,N}^{j,k}|\right]$ is $o(T_N^{1/\lambda})$
for $k\in [-2,2]$, $j\in \{0,\ldots,d\}$. 
Let us remind that we are interested in $k\in\{-2\gamma_j,1-2\gamma_j,2-2\gamma_j,2\gamma_j\}\cup \{-1,0,1,2\}$.

Let us first note that
\begin{equation}
    \calQ_{T_N}^{j,k}-\calQ_{T_N,N}^{j,k}
    =\int_0^{T_N}(X_t^k-X^k_{\lfloor t\rfloor_{\Delta_N}})\mathds{1}_{I_j}(X_t)\vd t +\int_0^{T_N}X^k_{\lfloor t\rfloor_{\Delta_N}}(\mathds{1}_{I_j}(X_t)-\mathds{1}_{I_j}(X_{\lfloor t\rfloor_{\Delta_N}}))\vd t
\end{equation}
and
\begin{equation}
    X^k_{\lfloor t\rfloor_{\Delta_N}}(\mathds{1}_{I_j}(X_t)-\mathds{1}_{I_j}(X_{\lfloor t\rfloor_{\Delta_N}}))
    = X_{\lfloor t\rfloor_{\Delta_N}}^k(\mathds{1}_{X_t\in I_j, X_{\lfloor t\rfloor_{\Delta_N}}\notin I_j}-\mathds{1}_{X_t\notin I_j, X_{\lfloor t\rfloor_{\Delta_N}}\in I_j }).
\end{equation}
Hence, in the notation of the previous section, we have
\begin{align} \label{Q}
& \Emu\left[|\calQ_{T_N}^{j,k}-\calQ_{T_N,N}^{j,k}| \right]
\leq \int_0^{T_N} \mathbf{J}_{k,j}(t) + \sum_{i=0}^d \mathbf{Y}_{k,i}(t) \vd t.
\end{align}
The terms involving $\mathbf{Y}_{k,i}$, via Proposition~\ref{Hitting_time}, lead to the bound~\eqref{eq:gsta_bound}. 
While we need to consider more in details the case
\begin{align} \label{eq:term:similar}
   \int_0^{T_N} \mathbf{J}_{k,j}(t) \vd t =\int_0^{T_N}\Emu\left[ |X_t^k-X^k_{\lfloor t\rfloor_{\Delta_N}}| \mathds{1}_{I_j}(X_t)  \right]\vd t.
\end{align}
This term is trivial if $k=0$, it is bounded from above by~\eqref{eq:1_condition_bound} if $k=1$. 
In the other cases, i.e.~$k\in [-2,0)\cup(0,1)\cup(1,2]$, note that $k$ appearing in $\mathbf{J}_{k,j}$ always depends on $j$, i.e.~$k=k_j$. 
We show that $\mathbf{J}_{k_j,j} \leq C \sqrt{\Delta_N}$ for some positive constant $C$.

Note that if $\gamma_0=0$, $\mu$ admits no negative moments smaller or equal than $-1$. So in what follows about negative moments, unless $\gamma \equiv 0$, assume $\gamma_0\neq 0$. If $\gamma \equiv 0$, then only $k=0,1,2$ have to be considered and no negative moment is involved. 

The values of $k_j$ are $2 \gamma_j$, $2(1-\gamma_j)$, $1-2\gamma_j$, $-2\gamma_j$, and $-1$.
In what follows we denote $\kappa_{j,i}:=k_j-2(1-\gamma_i)$. 
Note that $\kappa_{j,i}\leq k_j$ for all $i$ and $\kappa_{j,i} \geq k_j-1$ if and only if $\gamma_i\geq 1/2$.

\paragraph{The case $k=0,1$.}
$\mu$ admits finite $k$-th moment. 

\paragraph{The case $k\in (0,1)\cup(1,2]$.} 
If $j\notin\{ 0,d$\} $\mu$ admits finite $k_j$-th moment, if $j=0$ $\mu$ admits finite $\min(0,k_0-1)$-th moment,
if $j=d$ $\mu$ admits finite $2 k_d \mathds1_{(0,1)}(k_d) + 2(k_d+\gamma_d-1) \mathds1_{(1,2]}(k_d)$-th 
moment. Note that this is weaker than asking for $2(1+\gamma_d)$-th moment.  

We only need to discuss the moment conditions. 
Here $k_j=2\gamma_j, 2(1-\gamma_j)$ and, when $\gamma_j\in (0,1/2)$, also $k_j=(1-2\gamma_j)$. 
We consider separately $j=d$ and $j\neq d$. 

If $j=d$, H\"older's inequality ensures that
    \begin{align} \label{eq:term:similar:4}
   \mathbf{J}_{k_d,d}(t) 
   \leq \Emu\left[ |X_t^{k_d}-X^{k_d}_{\lfloor t\rfloor_{\Delta_N}}| \mathds{1}_{I_d}((X_u)_{u\in[\lfloor t\rfloor_{\Delta_N},t]})  \right] 
   + 2 \sum_{i=0}^d \Emu\left[ |X_t^{k_d}|^2+|X^{k_d}_{\lfloor t\rfloor_{\Delta_N}}|^2 \right]^{1/2}  (\mathbf{Y}_{0,i})^{1/2}.
    \end{align} 
    If $k_d \in (1,2]$ (and also $k_d\in (0,1)$), Proposition~\ref{prop:Ito} provides the bound for the first summand of the right-hand side of the latter equation, if $\mu$ admits finite moments of order $k_d, k_d-1, \kappa_{d,d}$, and $\kappa_{d,d} + k_d $ on $I_d$ and of order $k_d, k_{d}-1, \kappa_{d,0}, \kappa_{d,0}+k_d$ on $I_0$. This leads to positive moment of order $k_d +\max(0,\kappa_{d,d})$. 
   If $k_d \in (0,1)$, Proposition~\ref{prop:Ito} provides the bound for the first summand of the right-hand side of the latter equation, if $\mu$ admits finite moments of order $k_d, k_d-1<0, \kappa_{d,d}\leq0$, and $\kappa_{d,d} + k_d\leq k_d$ on $I_d$.
   If $\mu$ admits finite moments of order $2 k_d$, Proposition~\ref{Hitting_time} ensures that $\mathbf{Y}_{0,i}\leq C \Delta_N$ for some constant $C\in (0,+\infty)$ independent of $N$. 

If $j\neq d$, then
\begin{align} \label{eq:term:similar:2}
   \mathbf{J}_{k_j,j}(t) \leq \Emu\left[ |X_t^{k_j}-X^{k_j}_{\lfloor t\rfloor_{\Delta_N}}| \mathds{1}_{I_j}((X_u)_{u\in[\lfloor t\rfloor_{\Delta_N},t]})  \right] 
   + \sum_{i=0}^d \left( r_d^{k_j} \mathbf{Y}_{0,i} + \mathbf{Y}_{k_j ,i} \right).
\end{align}
Proposition~\ref{prop:Ito} gives the bound for the first term and Proposition~\ref{Hitting_time} yields the bound for the second term under the assumption that $\mu$ admits finite moments of order $k_j$ and if $j=0$ one should ad the moments of order $k_0, k_0-1, \kappa_{0,0}, \kappa_{0,0} +k_0$ on $I_0$. Note that $\kappa_{0,0}$ is always non-negative.

\paragraph{The case $k\in [-2, 0)$.} 

If $j\neq 0$, $\mu$ admits finite $k_j$-th moment. 
If $j=0$ and $k_0=-1$ then 
$\mu$ admits finite moments of order $-q_{\Lk}$ and $p_{\Lk}$ with $p_{\Lk},q_{\Lk}\geq 1$ such that $1=1/p_{\Lk} + 2/q_{\Lk}$. 
If $j=0$ and $k_0\in [-2,-1)\cup (-1,0)$, $\mu$ admits finite moments of order $-2(1+\gamma_0)$.  

Here $k_j=-1,-2\gamma_j$ and, when $\gamma_j\in (1/2,1]$, also $k_j=1-2\gamma_j$. 
We consider separately $j=0$ and $j\neq 0$. 

If $j\neq 0$, then
\begin{align} \label{eq:term:similar:3}
   \mathbf{J}_{k_j,j}(t) \leq \Emu\left[ |X_t^{k_j}-X^{k_j}_{\lfloor t\rfloor_{\Delta_N}}| \mathds{1}_{I_j}((X_u)_{u\in[\lfloor t\rfloor_{\Delta_N},t]})  \right] 
   + \sum_{i=0}^d \left( r_1^{k_j} \mathbf{Y}_{0,i} + \mathbf{Y}_{k_j ,i} \right).
\end{align}
Proposition~\ref{prop:Ito} gives the bound for the first term and Proposition~\ref{Hitting_time} yields the bound for the second term under the assumption that $\mu$ admits finite moments of order $k_j$. 

If $j=0$, we distinguish between $k_0=-1$ and $k_0\in [-2,-1)\cup(-1,0)$. 
In this latter case, by Proposition~\ref{prop:Ito} we get the bounds under the assumptions that $\mu$ admits finite moment of order $\min(k_0-1, k_0+\kappa_{0,0})$-th moment on $I_0$. These moments are all greater than $- 2(1+\gamma_0)$. 
If $k_0=-1$  
H\"older's inequality yields
\begin{align}
    & \Emu\left[ |X_t^{k}-X^{k}_{\lfloor t\rfloor_{\Delta_N}}|  \mathds{1}_{I_0}(X_t)
    \right] 
    \\
    & \leq 
    \Emu\left[ |X_t^{-k}-X^{-k}{\lfloor t\rfloor_{\Delta_N}}|^{p_{\Lk}} \mathds{1}_{ I_0 } (X_t) \right]^{1/p_{\Lk}}
    \Emu\left[ |X_t|^{k q_{\Lk}}\right]^{1/q_{\Lk}}
    \Emu\left[ |X_{\lfloor t\rfloor_{\Delta_N}}|^{k q_{\Lk}}\right]^{1/q_{\Lk}}. 
\end{align}
We conclude with~\eqref{eq:1_condition_bound}.  

\subsubsection{End of the proof of Lemma~\ref{Convergence_Tps_long}}
\label{sec:lemma52_proof}
In this section, we complete the proof of Lemma~\ref{Convergence_Tps_long} by summarising the assumptions of the previous sections for each item. 
We have seen in the previous sections that the unique condition on $\Delta_N$ is $\lim_{N\to \infty} T_N^{\lambda-1} \Delta_N=0$.

\paragraph{Proof of~Item~\ref{item_1} (QMLE $\theta_{T_N,N}^{(\qLk)}$).}
The QMLE~\eqref{Estim_discret_multi_QMLE_1} involves the statistics $\calQ^{j,k}$ and $\calM^{j,m}$, where $m \in \{0,1\}$, $k \in \{0,1,2\}$, and $j \in \{0,\ldots,d\}$. 
The most restrictive assumptions on the moments of $\mu$ is obtained for $m=1$: $\mu$ has finite moment of order $2(1+\gamma_d)$, that is $\calH_{\qLk}$.

\begin{remark} \label{rem:ass:QMLE}
   If $\sigma_\star$ is known, one could replace $\calM^{j,1}_{T_N,N}$ by $\frM^{j,1}_{T_N,N}$ in~\eqref{rewrite_Q_di} obtaining the analog of~\eqref{Estim_discret_multi_QMLE_1}. 
   However, in order to exploit the asymptotic properties of the estimator from continuous time observations, we need $\mathbf{H}_{\qLk}\subseteq \calH_{\qLk}$, hence $\calH_{\qLk}$ cannot be relaxed. 
  \end{remark}

\paragraph{Proof of~Item~\ref{item_0} (Volatility estimation).}  
Only the statistics $\calQ^{j,k}$, $\calM^{j,m}$, and $\calM^{i,0}$ with $m\in\{0,1\}$, $k=2\gamma_j$, and $i,j\in \{0,\ldots,d\}$ appear in the volatility estimator~\eqref{estim_sigma_discret}.
We observe that the most restrictive assumptions on the positive moments of $\mu$ is obtained for $m=1$: $\mu$ has finite moment of order $2(1+\gamma_d)$, that is $\calH_{\qLk}$.  

\paragraph{Proof of Item~\ref{item_2} (MLE $\theta_{T_N,N}^{(\Lk)}$ with $\sigma_\star$ unknown).}

To deal with the case $\sigma_\star$ is unknown, one takes the most restrictive conditions among the ones in Item~\ref{item_0} and the ones obtained in the case $\sigma_\star$ known.
First, assume that $\sigma_\star$ is known. 

The MLE~\eqref{Estim_discret_multi_MLE} involves the statistics $\calQ^{j,k}$, $\calM^{i,0}$, where $k \in \{-2\gamma_j, 1-2\gamma_j, 2-2\gamma_j\} \cup \{-1, 0\}$, and $i, j \in \{0, \ldots, d\}$. 
Then, $k_j=-2\gamma_j$ leads to the most restrictive negative moment conditions and $k_j=2(1-\gamma_j), -2\gamma_0$ to the positive ones.

If $\gamma_j=0$, then we get moments of order $0,1,2$, indeed $k_j=-1$ would not appear. In this case, we ask for finite second moment for $\mu$. 

If $\gamma\not\equiv 0$, then we would at least get $k_j=-1$ for some $j=0,\ldots,d$ and so, by Proposition~\ref{Moment_sup}, $\gamma_0=0$ is excluded.

If $\gamma_0=1/2$, let $p_{\Lk}, q_{\Lk}$ such that $1/p_{\Lk}+2/q_{\Lk}=1$. 
$p_{\Lk}$ and $q_{\Lk}$ appear only when $\gamma_0=1/2$.
The most restrictive negative moment is $-q_{\Lk}<-2$ and the positive one is the maximum among $p_{\Lk}$ and $2(1-\gamma_j), j=0,\ldots,d$. 

If $\gamma_0>1/2$, we are only interested in the positive moments conditions because Proposition~\ref{Moment_sup} ensures the process admits negative moments of all order. 
The most restrictive positive moment $\max_{j=0,\ldots,d}2(1-\gamma_j)$.  

We are now ready to combine with the conditions coming from the fact that $\sigma_{\star}$ is unknown.
The positive moments: $\max_{j}2(1-\gamma_j) \leq 2\leq 2(1+\gamma_d)$ but we cannot say for $p_{\Lk}$, it depends on $q_{\Lk}$.

\subsection{Proof of Proposition \ref{Preuve_holder}} \label{proof:Preuve_holder}
 Given $s$ and $t$ such that $0\leq s<t<s+1$, we show that for every $m\geq 1$ such that $\mu$ admits finite $m$-th moment, there exists a constant $C\in(0,\infty)$ depending only on $m$ and the parameters of the process such that $\Emu\left[|X_{t}-X_s|^m\right]\leq C(t-s)^{\sfrac{m}{2}}$.

By the triangular inequality,
\begin{align}
|X_t-X_s|&\leq \int_s^t\left|a(X_u)-b(X_u)X_u\right|\vd u+\left|\int_s^t\sigma(X_u)(X_u)^{\gamma(X_u)}\vd B_u\right|\\
&\leq (t-s)\max_{i=0,\ldots,d} |a_i|+\max_{i=0,\ldots,d} |b_i|\int_s^t|X_u|\vd u+\left|\int_s^t\sigma(X_u)(X_u)^{\gamma(X_u)}\vd B_u\right|.
\end{align}
Then, Jensen's inequality ensures that for $m\geq 1$ it holds that
\begin{align}
\Emu\left[|X_t-X_s|^m\right]&\leq 2^{2m-2}\max_{i=0,\ldots,d} |a_i|(t-s)^m+2^{2m-2}\max_{i=0,\ldots,d} |b_i|(t-s)^{m-1}\int_s^t\Emu\left[|X_u|^m\right]\vd u\\
&+2^{m-1}  \Emu\left[\left|\int_s^t\sigma(X_u)(X_u)^{\gamma(X_u)}\vd B_u\right|^m\right].
\end{align}

Since $X_0$ is distributed as the stationary distribution $\mu$, which admits finite $m$-th moment, then $\sup_{u\in[s,t]}\Emu\left[|X_u|^m\mathds{1}_{I_d}(X_u)\right]<\infty$.

Burkholder-Davis-Gundy inequality implies that
\[
    \mathcal{J}(m)
    :=\Emu\left[\left|\int_s^t\sigma(X_u)(X_u)^{\gamma(X_u)}\vd B_u\right|^m\right]\leq \Emu\left[\left(\int_s^t\max_{i=0,\ldots,d} |\sigma_i| (X_u)^{2\gamma(X_u)}\vd u\right)^{\!\!\sfrac{m}{2}}\right].
\] 
If $m\geq 2$, H\"older's inequality together with the fact that $\sup_{u\in[s,t]}\Emu\left[|X_u|^{m\gamma_d}\mathds{1}_{I_d}(X_u)\right]<\infty$,
yields 
\begin{align}
      \mathcal{J}(m) 
      &\leq (t-s)^{\sfrac{m}{2}-1}(\max_{i =0,\ldots,d}|\sigma_i|)^{\sfrac{m}{2}}\int_s^t\Emu\left[|X_u|^{m\gamma_d} \right]\vd u\nonumber 
    \leq C(t-s)^{\sfrac{m}{2}}.
\end{align} 
If $m\in[1,2)$, by H\"older's inequality
\begin{equation}
    \mathcal{J}(m) \leq \sqrt{  \mathcal{J}(2m) }
   \leq C (t-s)^{\sfrac{m}{2}}.
\end{equation}
The proof is thus completed.

\subsection{Proof of Proposition~\ref{prop:Ito}}

The techniques are similar to those in Section~\ref{proof:Preuve_holder}. 
It\^o-formula ensure that there exists constants $C_1,C_2 \in (0,+\infty)$ depending only on $a,b,\sigma$ such that
\begin{align}
     & \Emu\left[|X_t^k-X_s^k| \mathds{1}_{I_j}(X_t) \right]
     \leq 
    C_1 |k| \int_s^t \Emu\left[|X_u|^k + |X_u|^{k-1} \right] \vd u  
    \\& + C_2 |k||k-1| \int_s^t \Emu\left[|X_u|^{k-2+2\gamma(X_u)} \right] \vd u 
    + 
    \Emu\left[ |k|  \left|\int_s^t X_u^{k-1} \sigma(X_u) X_u^{\gamma(X_u)} \vd B_u \right| \right]
\end{align}
Burkholder-Davis-Gundy inequality ensures that there exists a constant $C_3\in (0,+\infty)$ such that
\begin{align}
        \Emu\left[ |k|  \left|\int_s^t X_u^{k-1} \sigma(X_u) X_u^{\gamma(X_u)} \vd B_u \right| \right]^2 
        \leq C_3 \int_s^t \Emu\left[ X_u^{2(k-1+\gamma(X_u))} \vd u \right].
\end{align}
Note that $I_i$ for all $i\in \{1,\ldots, d-1\}$ is contained in $[r_1,r_d] \subseteq (0,+\infty)$,
hence there exists a constant $C\in (0,+\infty)$ depending on $k,a,b,\sigma$ such that
\begin{align}
    & \Emu\left[|X_t^k-X_s^k| \mathds{1}_{I_j}(X_t) \right]
     \leq 
    C  (t-s) + C \int_s^t \left(\nu(k,d) + \nu(k-1,d) + \nu(k-2(1-\gamma_d),d) \right)  \vd u
    \\
    & \quad + C \int_s^t \left(\nu(k,0) + \nu(k-1,0) + \nu(k-2(1-\gamma_0),0)  \right)  \vd u
    \\
    & \qquad + C \sqrt{\int_s^t \nu(2(k-1+\gamma_d), d) + \nu(2(k-1+\gamma_0), 0) \vd u}.
\end{align}
where we denoted by $\nu(m,j)=\Emu\left[ X_u^{m} \mathds{1}_{I_j}(X_u)\right] = \int_{I_j} x^m \vd \mu(x)$. This proves the first statement. 
The proof of the second statement is similar.

\subsection{Proof of the key result: Proposition \ref{Hitting_time}} \label{proof:Hitting_time}

For all $j\in\{0,\ldots,d\}$, let $\xi^{(j)}$ denote a standard CKLS process with parameters $(a_j,b_j,\sigma_j,\gamma_j)$ starting at $X_s$. Let $s,t\in[0,\infty)$ be fixed such that $0\leq s<t <s+1$. Note that:
\begin{align}
    & \Emu\left[|X_s|^m \mathbb{P}_{X_s}\left(\tau_{I_j}^{\xi^{(j)}}<t-s\right)\right]
    \\
    &\leq  \Emu\left[|X_s|^m \left(\mathbb{P}_{X_s}\left( \tau_{r_j,\searrow}^{\xi^{(j)}}<t-s\right) +\mathbb{P}_{X_s}\left( \tau_{r_{j+1},\nearrow}^{\xi^{(j)}}<t-s\right)  \right)\mathds{1}_{I_j}(X_s)\right] ,
\end{align}
where $\tau_{r_j,\searrow}^{\xi^{(j)}}$ is the first hitting time from above of the level $r_j$, $\tau_{r_{j+1},\nearrow}^{\xi^{(j)}}$ is the first hitting time from below of the level $r_{j+1}$ of the process $\xi^{(j)}$.

Let $k\geq 1$, without loss of generality, we reduce to show that: 
\begin{equation}\label{tps_att_inf}
\Emu\left[|X_s|^m\mathbb{P}_{X_s}\left(\tau_{r_{d},\searrow}^{\xi^{(d)}}<t-s\right)\mathds{1}_{I_{d}}(X_s)\right]\leq C_1 (t-s)^{\sfrac{k}{2}},
\end{equation}
and
\begin{equation}\label{tps_att_sup}
\Emu\left[|X_s|^m\mathbb{P}_{X_s}\left(\tau_{r_{1},\nearrow}^{\xi^{(0)}}<t-s\right)\mathds{1}_{I_{0}}(X_s)\right]\leq C_2 (t-s)^{\sfrac{k}{2}},
\end{equation}
where $C_1$ and $C_2$ are strictly positive constant.

Indeed, in the other cases, $X_s$ belongs in $I_j$ for $j\notin \{0,d\}$, which is compact, and the desired inequality can be deduced using a similar reasoning.

\subsubsection{Bounds on the first hitting time from above of the level $r_d$: \eqref{tps_att_inf}}

We focus on the case $\gamma_d\in[0,1)$. 
The case $\gamma_d=1$ can be proven using a similar reasoning, the proof is thus omitted.
Let us recall that the parameters $(a_d,b_d,\sigma_d,\gamma_d)$ satisfy the ergodicity conditions in Table~\ref{table_ergo} in Section~\ref{sec:property_stationary}, in particular $b_d \geq 0$.\\
The main idea of this proof is to bound the first hitting time by the hitting times of some drifted Brownian motions. To do so, we apply the Lamperti transform and we bound the process, over a well chosen time interval.

We define the process $(Y_u)_{u\geq 0}$ as follows. For all $u\geq 0$ let $Y_u=\psi(\xi^{(d)}_u)$ where \\$\psi(x)=\int_{0}^x \frac{1}{\sigma_d y^{\gamma_d}}\vd y=\frac{x^{1-\gamma_d}}{\sigma_d(1-\gamma_d)}$ (Lamperti transform). We denote $\psi^{-1}$ as the inverse function of $\psi$, then $Y$ is solution to the following SDE:
\[
    \vd Y_u=\frac{a_d}{\sigma_d} (\psi^{-1}(Y_u))^{-\gamma_d}-b_d(1-\gamma_d)Y_u
    -\frac12 \sigma_d \gamma_d (\psi^{-1}(Y_u))^{\gamma_d-1}\vd u+\vd B_u.
\]

Let $\varepsilon >0$ be fixed and let $\tau^{Y}_{[\psi(r_d),\psi(X_s)+\varepsilon]}$ denote the first hitting time of the boundary of $[\psi(r_d),\psi(X_s)+\varepsilon]$ of the process $Y$.
Then, we have  
\[
    \mathbb{P}_{X_s}\left(\tau_{r_d,\searrow}^{\xi^{(d)}}<t-s\right)=\mathbb{P}_{\psi(X_s)}\left(\tau_{\psi(r_d),\searrow}^Y<t-s \right) 
    \leq \mathbb{P}_{\psi(X_s)}\left(\tau_{[\psi(r_d),\psi(X_s)+\varepsilon]}^Y<t-s \right).
\]

By applying the Comparison Theorem (e.g.,~\cite[p352]{ikeda2014stochastic}) until time 
$\tau^{Y}_{[\psi(r_d),\psi(X_s)+\varepsilon]}$, $Y-\psi(X_s)$ is bounded from above by a drifted BM $B^{\nu_+}$ starting from 0 and from below by a drifted BM $B^{\nu_-}$ starting from 0 with parameters:
\[
    \left\{
    \begin{array}{ll}
        \nu_+=\frac{|a_d|}{\sigma_d}{r_d^{-\gamma_d}} \\
        \nu_-=-\frac{|a_d|}{\sigma_d}{r_d^{-\gamma_d}} - b_d(1-\gamma_d)(\psi(X_s)+ \varepsilon)-\frac12 \sigma_d \gamma_d {r_d^{\gamma_d-1}}.
    \end{array}\right. 
\]
 Hence, the following inequality holds:
\[
    \mathbb{P}_{X_s}\left(\tau_{r_d,\searrow}^{\xi^{(d)}}<t-s\right)\leq  \mathbb{P}_{0}\left(\tau_{\psi(r_d)-\psi(X_s),\searrow}^{B^{\nu_-}}<t-s\right)+ \mathbb{P}_{0}\left(\tau_{\varepsilon,\nearrow}^{B^{\nu_+}}<t-s\right),
\]
and classical results on the first hitting of a drifted Brownian motion (see \cite{borodin2015handbook}) yield 
\begin{equation}
 \mathbb{P}_{0}\left(\tau_{\varepsilon,\nearrow}^{B^{\nu_+}}<t-s\right)\leq K_1e^{-\frac{\varepsilon^2}{2(t-s)}},
\end{equation}
and,
\begin{equation}
\mathbb{P}_{0}\left(\tau_{\psi(r_d)-\psi(X_s),\searrow}^{B^{\nu_-}}<t-s\right)\leq K_2e^{-\frac{(\psi(r_d)-\psi(X_s))^2}{2(t-s)}+b_d(1-\gamma_d)\psi(X_s)^2} f(X_s),
\end{equation}
where $K_1$ and $K_2$ are two strictly positive constants, $f$ is an explicit function which depends only on $X_s$ and such that $\lim_{x\to +\infty}f(x) e^{-\psi(x)^2} =0$. 

Let us note that
\[
    \mu(x)\mathds{1}_{x\geq r_d}=K_3\frac{2}{\sigma_d^2x^{1-2\frac{a_d}{\sigma_d^2}}}\exp\left(-b_d(1-\gamma_d)\psi(x)^2\right)\mathds{1}_{x\geq r_d},
\]
with $K_3$ a strictly positive constant.

Then, as $ s \to t $, and since $ \mu $ has a finite $ m $-th moment by assumption, we derive the following bound:
\begin{align}
    &\mathbb{E}\left[|X_s|^m \mathbb{P}_{X_s}\left(\tau_{r_{d},\searrow}^{\xi^{(d)}} < t - s\right) \mathds{1}_{I_{d}}(X_s)\right]  \\
   & \leq C\left(\int_{r_d}^{+\infty} x^m e^{-\frac{(\psi(r_d)-\psi(x))^2}{2(t-s)}} \vd x + e^{-\frac{\epsilon^2}{2(t-s)}}\int_{r_d}^{+\infty} x^m \mu(x) \vd x\right) \\  
    &\leq C\left( \int_{r_d}^{r_d+1} e^{-\frac{(\psi(r_d)-\psi(x))^2}{2(t-s)}} \vd x + \int_{r_d+1}^{+\infty} x^m e^{-\frac{(\psi(r_d)-\psi(x))^2}{2(t-s)}} \vd x + e^{-\frac{\epsilon^2}{2(t-s)}} \right) \\  
    &\leq C \left(e^{-\frac{(\psi(r_d) - \psi(r_d+1))^2}{2(t-s)}} + \int_{\frac{\psi(r_d+1) - \psi(r_d)}{\sqrt{t-s}}}^{+\infty} \psi^{-1}\left(x\sqrt{t-s} - \psi(r_d)\right)^m e^{-x^2} \vd x + e^{-\frac{\epsilon^2}{2(t-s)}} \right) \\  
    &\leq C(t-s)^{\frac{k}{2}}.  
\end{align}
where $ C $ is a strictly positive generic constant. The last inequality follows from an asymptotic expansion of the tail probabilities of a Gaussian distribution. This completes the proof. 

Note that, in the preceding calculations, when $ (t - s) $ is close to zero, the term $ e^{-K/(t-s)} $ decays faster than any polynomial of $ t - s $. To handle this, we split the integral to distinguish between cases where $ \psi(x) - \psi(r_d) $ is small (e.g., of order $ O(\sqrt{t-s}) $) and cases where it is large (e.g., of order $ O(1/\sqrt{t-s}) $). In the latter case, we can always bound the integral using the asymptotic expression of the tail probabilities of a Gaussian distribution

From now on, we will systematically apply this approach and omit the detailed calculations.

\subsubsection{Bounds on the first hitting time from below of the level $r_1$: \eqref{tps_att_sup}} \label{sec:proof:around0}
The parameters $(a_0,b_0,\sigma_0,\gamma_0)$ satisfy the conditions ensuring ergodicity in Table~\ref{table_ergo}.
 We remark that, on a suitable time interval, the process $(\xi^{(0)})^{2(1-\gamma_0)}$ can be bounded from above by the norm of a multi-dimensional Brownian motion. This leads to obtaining \eqref{tps_att_sup} under the assumption that $\gamma_0\in \{0,{1}/{2}\}$. Instead, if $\gamma_0\in ({1}/{2},1]$, this bound is not enough. Thus, we additionally bound from below the Lamperti transform of the process $\xi^{(0)}$ by a drifted Brownian motion.

\paragraph{The case $\gamma_0={1}/{2}$.}

We remind the following results, for $i\in\{1,\ldots,n\}$, we denote $\tilde{B}:=(\tilde{B}^i)_{i\leq n}$ a n-dimensional Brownian motion. It\^o formula and Levy characterization imply that 
$$ \forall u\geq 0,\quad \left|\left| \frac{\sigma_0}{2}\tilde{B}_u-\sqrt{\frac{X_s}{n}}\right|\right|^2_2=Y_{u,n}=X_s+n\frac{\sigma_0^2}{4}u+\int_0^u\sigma_0\sqrt{Y_{v,n}}\vd W_v,$$
where $W$ is a Brownian motion.
Then $Y_{\cdot,n}$ is a CIR process whose coefficients satisfy the conditions for ergodicity in Table~\ref{table_ergo}. 
Moreover we take $n$ such that $0<a_0<{n\sigma_0^2}/{4}$. So, by the Comparison Theorem, it holds a.s.~that $\xi_u^{(0)}\leq Y_{u,n}$ for all $0\leq u\leq \tau_{r_1,\nearrow}^{\xi^{(0)}}$ and then
\begin{equation}
    \mathbb{P}_{X_s}\left(\tau_{r_1,\nearrow}^{\xi^{(0)}}\leq t-s\right)\leq \mathbb{P}_{X_s}\left(\tau_{r_1,\nearrow}^{Y_{.,n}}\leq t-s\right).
\end{equation}
Moreover, we have
\begin{equation}
    \{\tau_{r_1,\nearrow}^{Y.,n}\leq t-s\}\subseteq \bigcup_{i=1}^{n} \left\{\forall u\leq t-s, -\sqrt{\frac{r_1}{n}}<\left(\frac{\sigma_0}{2}\tilde{B}_u^i-\sqrt{\frac{X_s}{n}}\right) < \sqrt{\frac{r_1}{n}}\right\}^c .
\end{equation}
Then, by the Comparison Theorem and the symmetry of BM, we obtain:
\begin{equation}
\mathbb{P}_{X_s}\left( \tau_{r_1,\nearrow}^{\xi^{(0)}}\leq t-s\right)
\leq 2n\mathbb{P}_{0}\left(\tau_{2\frac{\sqrt{r_1}-\sqrt{X_s}}{\sigma_0\sqrt{n}},\nearrow}^{\tilde{B}^1}\leq t-s\right)\leq K_4e^{-\frac{2\left(\sqrt{r_1}-\sqrt{X_s}\right)^2}{\sigma_0^2 n(t-s)}},
\label{tps_brownien}
\end{equation}
with $K_4$ a strictly positive constant.
This, and the stationary distribution $\mu$~\eqref{mesure_stationnaire}, if $\mu$ admits finite $m$-th moment, yield inequality~\eqref{tps_att_sup}.

\paragraph{The case $\gamma_0=0$.}
By applying It\^o formula on $Y:=(\xi^{(0)})^{2}$, for all $u\geq 0$, we have:
\begin{equation}
    \vd Y_u
    = 2\left(\operatorname{sgn}(\xi^{(0)}) a_0\sqrt{Y_u}-b_0Y_u+\frac{\sigma_0^2}{2}\right)\vd u+2\sigma_0\sqrt{Y_u} \vd \tilde B_u
\end{equation}
with $Y_0=X_s^{2}$ and $\tilde{B}_u$ another Brownian motion.
Since $2\operatorname{sgn}(\xi^{(0)}) a_0\sqrt{Y_u}-2b_0Y_u+\sigma_0^2 \leq 2|a_0| r_1 + \sigma_0^2 \leq n \sigma_0^2$ for some $n\in \mathbb{N}$, the Comparison Theorem ensures that for all $u\leq \tau^{Y}_{r_1^{2},\nearrow}$ the process $Y$ is bounded from above by the norm of a $n$-dimensional Brownian motion. Hence, similarly to the case $\gamma_0=\frac{1}{2}$, we have:
\begin{equation} \label{tps_brownien2}
    \mathbb{P}_{X_s}\left( \tau_{r_1,\nearrow}^{\xi^{(0)}}\leq t-s\right)\leq 2n\mathbb{P}_{0}\left(\tau_{\frac{r_1-X_s}{\sigma_0\sqrt{n}},\nearrow}^{\tilde{B}^1}\leq t-s\right)\leq K_5e^{-\frac{\left(r_1-X_s\right)^2}{2\sigma_0^2n(t-s)}},
\end{equation}
with $K_5$ a strictly positive constant. We conclude analogously to the case 
$\gamma_0={1}/{2}$.

\paragraph{The case $\gamma_0\in (1/2,1)$.}
By Lamperti transform, $Y_u :=\vd \psi(\xi_u^{(0)})$ with $\psi(x)=\frac{x^{1-\gamma_0}}{\sigma_0(1-\gamma_0)}$ (in particular $Y_0=\frac{X_s^{1-\gamma_0}}{\sigma_0(1-\gamma_0)}$). So,
\begin{align}
    \vd Y_u =\left[\frac{a_0}{\sigma_0}\left((1-\gamma_0)Y_u\right)^{1-\frac{1}{(1-\gamma_0)}}-b_0(1-\gamma_0)Y_u-\frac{\gamma_0}{2(1-\gamma_0)}Y_u^{-1}\right]\vd u+\vd B_u.
\end{align}
By the Comparison Theorem, it holds a.s.~for all $u\in[0, \tau_{\left[\frac{\psi(X_s)}{2},\psi(r_1)\right]})$ that
$$B_u^{\nu_-}\leq Y_u  \quad \text{and} \quad Y_u^2\leq \left|\left| \frac{1}{2}\tilde{B}_u-\frac{\psi(X_s)}{\sqrt{n}}\right|\right|^2_2,$$
where $B^{\nu_-}$ is drifted Brownian motion and $\tilde{B}:=(\tilde{B}^i)_{i\leq n}$ a n-dimensional Brownian motion. Here the drift parameter $\nu_-$ and the dimension $n$ are  given by
$$\left\{
    \begin{array}{ll}
        \nu_-=-\frac{|b_0|r_1^{1-\gamma_0}}{\sigma_0}-\frac{\gamma_0\sigma_0}{X_s^{1-\gamma_0}}, \\
        n(X_s)=\left\lceil \frac{a_0}{\sigma_0^2(1-\gamma_0)}X_s^{1-2\gamma_0}+|b_0|(1-\gamma_0)\frac{r_1^{2(1-\gamma_0)}}{\sigma_0^2(1-\gamma_0)^2} \right\rceil .
    \end{array}\right. $$

Finally, we obtain the following inequality: 
\begin{align}
\mathbb{P}_{X_s}\left(\tau_{r_1,\nearrow}^{\xi^{(0)}}<t-s\right)&\leq  \mathbb{P}_{0}\left(\tau_{-\frac{\psi(X_s)}{2},\searrow}^{B^{\nu_-}}<t-s\right)+ 2n(X_s)\mathbb{P}_{0}\left(\tau_{2\frac{\psi(r_1)^2-\psi(X_s)^2}{\sqrt{n(X_s)}},\nearrow}^{\tilde{B}^1}<t-s\right)\nonumber \\
&\leq K_6\left( e^{-\frac{|\psi(r_1)-\psi(X_s)|^2}{2(t-s)}}+n(X_s)e^{-\frac{2\left(\psi(r_1)^2-\psi(X_s)^2\right)^2}{n(X_s)(t-s)}}\right),
\label{tps_brownien3}
\end{align}
with $K_6$ a strictly positive constant.

Similarly to the previous cases, this, and the stationary distribution $\mu$~\eqref{mesure_stationnaire}, if $\mu$ admits finite $m$-th moment, yield inequality~\eqref{tps_att_sup}.

\paragraph{The case $\gamma_0=1$.}
One can conclude by bounding directly the Lamperti transform from above by a drifted Brownian motion as it has been done for the bounds on the first hitting time from above of the level $r_d$.

\paragraph{Acknowledgements.} The authors are grateful to P.~Pigato for discussions and suggestions about this work. They also thank A.~Lejay, and L.~Lenôtre for references on existing bounds for first hitting times and C.~Blanchet-Scalliet and D.~Dorobantu for discussions on related topics.
The authors thank the referee for his helpful remarks and suggestions.

\bibliographystyle{abbrv}
\bibliography{articleR1_arxiv_bib}

\end{document}